\pdfoutput=1
\documentclass[11pt,letterpaper]{article}

\usepackage[margin=1in]{geometry}
\usepackage{amsmath,amssymb,amsthm}
\usepackage{mathpazo}
\usepackage{graphicx}
\usepackage{booktabs}
\usepackage{enumitem}
\usepackage{hyperref}
\usepackage{natbib}
\usepackage{needspace}
\usepackage{caption}
\usepackage{xcolor}
\usepackage{tikz}
\usetikzlibrary{positioning,calc}
\usepackage{tikz-cd}

\hypersetup{colorlinks=true,linkcolor=black,citecolor=black,urlcolor=blue}

\newtheorem{remark}{Remark}

\title{Decorated Cospans at Work:\\Coupling Heterogeneous Dynamical Systems\\via Pushouts and Particle Filters}
\author{Wesley Phoa\thanks{Corresponding author: \texttt{wesley@topos.institute}.
Claude Opus 4.6 (Anthropic) was used for exploration, implementation, and editorial assistance.}}
\date{March 2026}

\begin{document}
\maketitle

\begin{abstract}
Decorated cospans provide a categorical framework for composing open systems along shared interfaces. This paper is a computational proof of concept: we show that the framework produces a working coupled dynamical system when the decorations are quantitative models from different mathematical traditions. Specifically, we couple a linearised New Keynesian DSGE, a stochastic compartmental epidemic (multi-strain SEIR), and a nonlinear vaccine adoption model with hysteresis into a single sequential Monte Carlo sampler. Each model is a decorated cospan---interior dynamics as decoration, exposed variables as interfaces. The composite system is the pushout along variable identifications, with coupling functions encoded as factor graph constraints. The pushout does not care what is inside the models; it cares only about what they share. The calibration is best-guess; the models deliberately simple; the aim is to show that the categorical machinery \emph{computes}, not to produce a forecasting tool.

The coupled system produces a rejection bifurcation: some trajectories escape via vaccination, others enter a self-reinforcing cycle of mandate backlash, vaccine refusal, sustained infection, and recession. This is a structural property of the coupling, not an input assumption. Coupling shifts the output gap by $-0.78$~pp and rejection by $+22$~pp relative to the uncoupled system. A fourth narrative---fiscal/political dynamics, calibrated to the US COVID fiscal response---attaches via a second pushout and introduces the first positive coupling channel. With pandemic-scale spending parameters, 14\% of trajectories overshoot into positive output gap territory; the bearish bias shrinks from $-0.78$~pp to $-0.20$~pp but persists. A computable bias decomposition separates three sources of this asymmetry---sampling, structural, and observational---and localises the structural component to specific coupling functions whose directional asymmetry can be tested against historical analogues. If realistically calibrated, one could formulate a response---a `game plan'---for each archetypal trajectory; as the actual historical path evolves, monitor its distance to each archetype across all state variables, and when it converges on one, implement the corresponding plan.
\end{abstract}

\medskip
\noindent\textbf{Keywords:} decorated cospans, pushouts, applied category theory, particle filters, coupled dynamical systems, factor graphs, DSGE, epidemiology, sequential Monte Carlo

\bigskip

\needspace{4\baselineskip}
\section{Introduction}

Decorated cospans~\citep{fong2016} offer an elegant algebra for composing open systems: each system exposes interfaces, and the pushout glues systems along shared interface variables while the interior dynamics ride as decorations, untouched by the composition. The theory is clean. But does it compute? When the decorations are not finite automata or Petri nets but quantitative dynamical systems---ODEs, stochastic processes, rational-expectations models---does the pushout produce a coupled system one can actually simulate?

This paper answers that question. The companion paper~\citep{phoa2026} develops a qualitative framework using decorated cospans, factor graphs, and particle filters for scenario analysis. Here we instantiate it computationally: three models from different mathematical traditions---a linearised New Keynesian DSGE, a stochastic compartmental epidemic, and a nonlinear vaccine adoption model with hysteresis---coupled through six factor-graph constraints into a single joint sequential Monte Carlo sampler. The models share no common mathematical structure. What they share is \emph{interfaces}: variables exposed for coupling. The pushout does not care what is inside.

The categorical point is that heterogeneous composition is the natural case, not the special case. A linearised rational-expectations model, a nonlinear ODE system, and a stochastic jump process compose using the same algebraic operation. The practical point is that the coupled system produces dynamics none of the components predicts alone: a rejection bifurcation, a compositional bearish bias of $-0.78$~pp, and structural asymmetries traceable to specific coupling functions. If realistically calibrated, one could formulate a game plan for each archetypal trajectory, with cross-narrative implications already traced by the categorical structure.

\needspace{4\baselineskip}
\section{Framework}\label{sec:framework}

We summarise the relevant machinery from~\citet{phoa2026}, adding only what the quantitative setting requires.

\needspace{3\baselineskip}
\subsection{Particle filters for single narratives}

A \emph{narrative} is a hidden Markov model $({\cal S}, {\cal Y}, f, g, \pi_0)$. A \emph{particle filter} approximates the filtering distribution by propagating $N$ weighted samples (not to be confused with the population normalisation $N$ in the SEIR model of Section~\ref{sec:seir}; context disambiguates):
\begin{align}
\text{Propagate:} \quad & \tilde{x}_t^{(i)} \sim f(\cdot \mid x_{t-1}^{(i)}) \\
\text{Reweight:} \quad & w_t^{(i)} \propto w_{t-1}^{(i)} \cdot g(y_t \mid \tilde{x}_t^{(i)}) \\
\text{Resample when:} \quad & \text{ESS} = \frac{1}{\sum_i (w_t^{(i)})^2} < \tfrac{N}{2} \quad \text{(where } \textstyle\sum_i w_t^{(i)} = 1\text{)}
\end{align}
Here the transition kernel \emph{is} the model: NK policy functions, SEIR differential equations, vaccine adoption dynamics.

\needspace{3\baselineskip}
\subsection{Factor graphs for multiple narratives}

When $k$ narratives interact, their dependencies are encoded in a factor graph~\citep{kschischang2001}: variable nodes (state-space variables) and factor nodes (coupling constraints). Each factor either modifies states (hard coupling, $f_j: {\cal S}_{i_1} \times {\cal S}_{i_2} \to {\cal S}_{i_1}$ or ${\cal S}_{i_2}$) or returns a compatibility score (soft coupling, $f_j: {\cal S}_{i_1} \times {\cal S}_{i_2} \to \mathbb{R}_{\ge 0}$).

\needspace{3\baselineskip}
\subsection{Decorated cospans for heterogeneous composition}

Each narrative is a decorated cospan~\citep{fong2016,fong2019} $(A_i \to N_i \leftarrow B_i, \delta_i)$:
\begin{itemize}[nosep]
\item $N_i$ is the interior state space
\item $A_i, B_i$ are interfaces---variables exposed for coupling
\item $\delta_i$ is the decoration: the dynamics (transition kernel, observation function)
\end{itemize}

The composite system is the pushout $N = \bigsqcup_i N_i / {\sim_{\cal V}}$ along variable identifications, with the product decoration restricted by factor constraints. The categorical point: \emph{the pushout is agnostic to the decorations}. A linearised DSGE, a nonlinear ODE, and a stochastic jump process compose along their interfaces using the same algebraic operation. The interior dynamics ride as decorations on a structure that cares only about what is shared.

\begin{remark}
The product decoration restricted by factor constraints is a valid decorated cospan composition in the sense of~\citet{fong2016}. The decoration functor $F: (\mathbf{FinSet}, +) \to (\mathbf{Set}, \times)$ is defined by $F(N) = \{\text{transition kernels on } N\} \times \{\text{observation functions on } N\}$, with the monoidal structure given by the product. The restriction by factor constraints corresponds to composing the decorations via the monoidal structure of~$F$, as in~\citet[Theorem~4.1]{fong2016}. The full verification is provided in the companion paper~\citep{phoa2026}. Briefly: define $F : (\mathbf{FinSet}, +) \to (\mathbf{Set}, \times)$ by $F(N) = \mathrm{Kern}(N) \times \mathrm{Obs}(N)$, where $\mathrm{Kern}(N)$ is the set of Markov kernels on $N$ and $\mathrm{Obs}(N)$ is the set of observation functions on $N$. The lax monoidal structure is given by the product: $F(N_1 + N_2) \to F(N_1) \times F(N_2)$ sends a kernel on the disjoint union to its restrictions. Composition of decorated cospans then produces the product kernel restricted by factor constraints, as in Definition~2. The factor constraints enter as a further restriction on the decoration, compatible with the monoidal structure because each constraint is local (involving at most two component state spaces).
\end{remark}

The categorical framework provides three concrete benefits beyond notational convenience. First, \emph{modularity}: the coupling architecture (which variables are shared, which factors constrain them) is specified independently of the component models. Replacing the linear New Keynesian model with a nonlinear DSGE model requires changing only the decoration on one cospan, not the coupling structure. Second, \emph{associativity}: Proposition~1 guarantees that hierarchical decomposition of a many-narrative system produces the same composite regardless of the order of composition, which underwrites the scalability of the approach. Third, \emph{the factor--variable distinction}: the categorical framework forces a separation between variable identification (pushout) and factor coupling (decoration), which corrected an error in our initial informal analysis where these were conflated (see the companion paper for details).

\needspace{4\baselineskip}
\section{The Three Narratives}\label{sec:models}

\needspace{3\baselineskip}
\subsection{Economic cycle (N$_1$): New Keynesian DSGE}

The standard three-equation New Keynesian model, linearised around steady state:
\begin{align}
\pi_t &= \beta \, \mathbb{E}_t[\pi_{t+1}] + \kappa \, y_t + \varepsilon_t^s \label{eq:phillips} \\
y_t &= \mathbb{E}_t[y_{t+1}] - \sigma^{-1}(i_t - \mathbb{E}_t[\pi_{t+1}] - r_t^n) \label{eq:is} \\
i_t &= \phi_\pi \, \pi_t + \phi_y \, y_t + \varepsilon_t^m \label{eq:taylor}
\end{align}
State variables: output gap~$y_t$, inflation~$\pi_t$, nominal interest rate~$i_t$. The underlying production function is $Y_t = A_t L_t$, where $L_t$ is the labour force; in the linearised model, $L$ disappears into the supply shock~$\varepsilon_t^s$. Exogenous shocks: supply~$\varepsilon_t^s$ (AR(1), persistence $\rho_s = 0.9$ quarterly), monetary policy~$\varepsilon_t^m$ (i.i.d.), natural rate~$r_t^n$ (AR(1), persistence $\rho_r = 0.8$ quarterly). AR(1) persistence parameters are converted from quarterly to weekly frequency: $\rho_s^{w} = \rho_s^{1/13} \approx 0.992$, $\rho_r^{w} = \rho_r^{1/13} \approx 0.983$. Calibration follows~\citet{gali2015}: $\beta = 0.99$, $\kappa = 0.024$, $\phi_\pi = 1.5$, $\phi_y = 0.125$.

Solved analytically for the minimum state variable solution: policy functions $x_t = P \, s_{t-1} + Q \, \eta_t$, where $s$ collects AR(1) shock states and $\eta$ collects innovations. Once extracted, the rational expectations solver is discarded and the model steps forward by matrix multiplication---weekly, with shock variances scaled by $\sqrt{1/13}$ from quarterly.

\needspace{3\baselineskip}
\subsection{Epidemic (N$_2$): SEIR with strain emergence}\label{sec:seir}

A compartmental SEIR model~\citep{kermack1927} with stochastic strain replacement, calibrated to COVID-19:
\begin{align}
\dot{S} &= -\beta_{\text{eff}} \, S_{\text{eff}} \, I / N + \omega \, R \\
\dot{E} &= \beta_{\text{eff}} \, S_{\text{eff}} \, I / N - \sigma \, E \\
\dot{I} &= \sigma \, E - \gamma \, I \\
\dot{R} &= (1-\mu) \, \gamma \, I - \omega \, R \\
\dot{D} &= \mu \, \gamma \, I
\end{align}
State variables: $S$, $E$, $I$, $R$, $D$ (population fractions satisfying $S + E + I + R + D = 1$), where $S_{\text{eff}}$ is the effective susceptible population after accounting for vaccination (defined in Section~\ref{sec:coupling}, Table~\ref{tab:coupling}). Fixed parameters: incubation $\sigma = 1.41$/week, recovery $\gamma = 0.7$/week, waning $\omega = 1/52$/week. Behavioural dampening: $\beta_{\text{eff}} = \beta(1 - \alpha I)$, $\alpha = 5$. In simulation, $S_{\text{eff}}$ is clamped to $[0, S]$ to prevent unphysical negative values at high vaccination coverage.

\textbf{Strain emergence.} New dominant strains arrive as a Poisson process, rate $\lambda = 1/40$/week ($\approx$10~months, matching Alpha$\to$Delta$\to$Omicron). Each draws $R_0 \sim \text{Uniform}(1.5, 6.0)$, immune escape $\varepsilon \sim \text{Beta}(3, 3)$, IFR $\mu \sim \text{Beta}(2, 40)$. On arrival, fraction~$\varepsilon$ of recovered return to susceptible. Single dominant strain at a time.

\needspace{3\baselineskip}
\subsection{Vaccine dynamics (N$_3$): innovation, adoption, rejection}

Three state variables: efficacy~$v \in [0,1]$, uptake rate~$u \in [0,1]$, rejection fraction~$\rho \in [0,1]$.
\begin{align}
v_{t+1} &= v_t + \Delta v_{\text{innovation}} - \Delta v_{\text{drift}} \\
\rho_{t+1} &= \rho_t + \theta_{\text{up}} \cdot m_t \cdot (1-\rho_t) - \theta_{\text{down}} \cdot (1-m_t) \cdot \rho_t \\
u_{t+1} &= u_t + \theta_{\text{adopt}} \cdot (\bar{u}(I_t) - u_t)^+ - \theta_{\text{decay}} \cdot u_t, \quad u \leq 1-\rho
\end{align}
Efficacy jumps on vaccine innovation (Poisson, rate $\lambda_v = 1/26$/week) and degrades on strain escape. Uptake targets $\bar{u}(I) = 0.3 + 2I$, rising with outbreak severity, but is capped by rejection: only non-rejecters can adopt.

The key dynamic is \textbf{rejection hysteresis}. Mandate pressure activates only when $I > 2\%$. Below that: no mandates, no ratchet. Above it, rejection grows at $\theta_{\text{up}} = 0.005$/week and decays at $\theta_{\text{down}} = 0.003$/week---a $5{:}3$ asymmetry that permanently shifts some fraction from hesitant to resistant. Calibrated loosely on US COVID hesitancy: initial $\rho_0 = 0.15$ (a novel pathogen's clean slate), rising and never fully returning.

\needspace{4\baselineskip}
\section{Cross-Narrative Coupling}\label{sec:coupling}

\begin{figure}[htbp]
\centering
\begin{tikzpicture}[
  var/.style={circle, draw, minimum size=8mm, inner sep=1pt, font=\footnotesize},
  econ/.style={var, fill=blue!12, draw=blue!40},
  epi/.style={var, fill=red!12, draw=red!40},
  vacc/.style={var, fill=green!12, draw=green!40},
  factor/.style={rectangle, fill=black, minimum size=4mm, inner sep=0pt},
  edge/.style={draw=gray!70, thick},
  ident/.style={draw=black, thick, dashed},
  flabel/.style={font=\scriptsize\itshape, text=black!70},
  collabel/.style={font=\bfseries\small},
]
\def\ecox{0}
\def\epix{3.5}
\def\vacx{7}

\node[collabel, blue!70] at (\ecox, 4.5) {Economy};
\node[collabel, red!70] at (\epix, 4.5) {Epidemic};
\node[collabel, green!50!black] at (\vacx, 4.5) {Vaccine};

\node[econ] (y)  at (\ecox, 3.5) {$y$};
\node[econ] (pi) at (\ecox, 2.0) {$\pi$};
\node[econ] (i)  at (\ecox, 0.5) {$i$};

\node[epi] (S)  at (\epix, 3.5) {$S$};
\node[epi] (I)  at (\epix, 2.0) {$I$};
\node[epi] (D)  at (\epix, 0.5) {$D$};

\node[vacc] (v)   at (\vacx, 3.5) {$v$};
\node[vacc] (u)   at (\vacx, 2.0) {$u$};
\node[vacc] (rho) at (\vacx, 0.5) {$\rho$};

\node[factor] (f1) at (1.75, 2.9) {};
\node[factor] (f2) at (1.75, 1.5) {};
\node[factor] (f4) at (5.25, 2.75) {};
\node[factor] (f5) at (3.5, -0.5) {};
\node[factor] (f6) at (3.5, -1.2) {};

\node[flabel, below left=1pt] at (f1) {$f_1$};
\node[flabel, below left=1pt] at (f2) {$f_2$};
\node[flabel, above=2pt] at (f4) {$f_4$};
\node[flabel, above=2pt] at (f5) {$f_5$};
\node[flabel, below=2pt] at (f6) {$f_6$};

\draw[edge] (I) -- (f1);
\draw[edge, ->] (f1) -- (y);

\draw[edge] (I) -- (f2);
\draw[edge, ->] (f2) -- (pi);

\draw[edge] (v) -- (f4);
\draw[edge] (u) -- (f4);
\draw[edge] (rho) -- (f4);
\draw[edge, ->] (f4) -- (S);

\draw[edge] (y.south) -- (0, -0.5) -- (f5);
\draw[edge, ->] (f5) -- (7, -0.5) -- (rho.south);

\draw[edge] (i.south) -- (0, -1.2) -- (f6);
\draw[edge, ->] (f6) -- (7, -1.2) -- (v.south);

\draw[ident] (y.north) to[out=60,in=150] node[above, font=\scriptsize\itshape] {$N \equiv L$} (S.north);

\end{tikzpicture}
\caption{Factor graph for the coupled narrative system. Variable nodes are coloured by narrative; factor nodes (black squares) encode cross-narrative coupling. The dashed arc indicates the variable identification $N \equiv L$ (population~$\equiv$ labour force). Factors $f_1$ and $f_2$ transmit epidemic shocks to the economy; $f_4$ transmits vaccine protection back to the epidemic; $f_5$ and $f_6$ transmit economic conditions to vaccine dynamics. Factor $f_3$ (outbreak $\to$ uptake demand) is embedded in the vaccine model's target function and does not appear as a separate node.}\label{fig:factor-graph}
\end{figure}
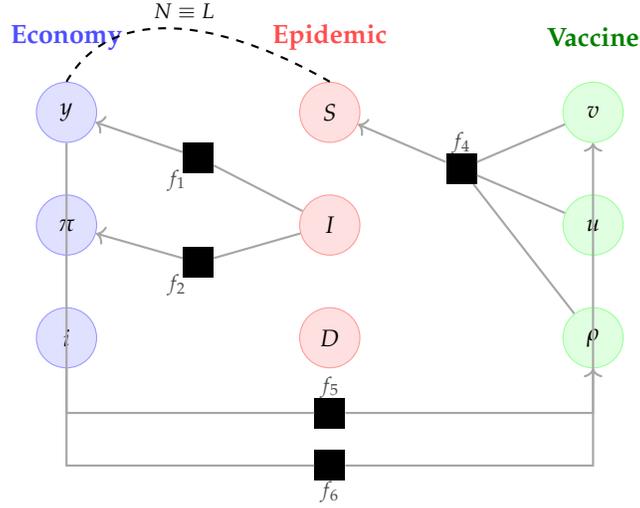

Six coupling functions connect the three narratives (Figure~\ref{fig:factor-graph}).

\begin{table}[htbp]
\centering
\caption{Cross-narrative coupling functions.}\label{tab:coupling}
\begin{tabular}{@{}llll@{}}
\toprule
Factor & Direction & Mechanism & Functional form \\
\midrule
$f_1$ & N$_2 \to$ N$_1$ & Demand shock & $\Delta r^n = -\eta_d(t) \cdot I_t$ \\
$f_2$ & N$_2 \to$ N$_1$ & Labour force & $\Delta \varepsilon^s = \eta_s(t) \cdot I_t$ \quad (absenteeism) \\
$f_3$ & N$_2 \to$ N$_3$ & Outbreak demand & Embedded in $\bar{u}(I)$ \\
$f_4$ & N$_3 \to$ N$_2$ & Effective immunity & $S_{\text{eff}} = S - v \cdot u \cdot (1-\rho)$ \\
$f_5$ & N$_1 \to$ N$_3$ & Mandate backlash & $\theta_{\text{up}} \times [1 + 5 \cdot \max(0, -y)]$ \\
$f_6$ & N$_1 \to$ N$_3$ & R\&D funding & $\lambda_v \times \max(0.5, 1 - 0.4 \cdot (i - i^*))$ \\
\bottomrule
\end{tabular}
\end{table}

\needspace{2\baselineskip}
\textbf{Calibration rationale.} The demand and supply channels ($f_1$, $f_2$) model the pure disease channel---not lockdowns. $f_1$ captures voluntary distancing and consumer fear; $f_2$ is the labour force channel---absenteeism, with $\xi = 0.3$ (30\% of infected workers absent at any given time). Both decay with \textbf{habituation}: $\eta_d(t) = 0.02 + 0.08\,e^{-0.02t}$, so the initial demand elasticity of 0.10 falls to 0.02 with a half-life of $\sim$35~weeks. Populations adapt to endemic disease: the first wave disrupts the economy; the fifth barely registers. Without this decay, the output gap never recovers---a calibration error that produced overly pessimistic conclusions in an earlier version of this analysis (see Discussion). Mandate backlash ($f_5$) triples rejection growth in a deep recession ($y = -10\%$). The R\&D channel ($f_6$) is deliberately mild.

\needspace{2\baselineskip}
\textbf{The mandate threshold.} Mandates fire only when $I > 2\%$---not at endemic levels. This creates \emph{intermittent} coupling: $f_5$ is silent during troughs and active during outbreaks. Combined with hysteresis, this produces path dependence: particles with severe early outbreaks trigger mandates, ratchet rejection upward, and never recover. Particles with mild early waves never trigger it.

\needspace{3\baselineskip}
\subsection{Variable identification: $N \equiv L$}

The SEIR population~$N$ and the NK labour force~$L$ are the same people. This is a pushout---two variables in different state spaces that are literally the same quantity. $f_2$ makes the identification operational: the effective labour force is $L_t = (1 - D_t) - \xi \, I_t$, where $D_t$ is cumulative mortality and $\xi$ the fraction of infected workers absent at any given time. The absenteeism component~$\xi I_t$ raises marginal costs---sick workers are a supply disruption, entering the Phillips curve as a cost-push shock. Deaths permanently reduce the production frontier but are not themselves cost-push. The demand channel~$f_1$ is separate: consumer fear and voluntary distancing enter the IS curve through the natural rate.

With an explicit labour market~\citep{gali2011}, the identification $N \equiv L$ would generate unemployment as a distinct state variable, and the decomposition into voluntary non-participation (caring for sick family), involuntary unemployment (firms laying off), and permanent exit (death) would matter for welfare. Worth doing---but the proof of concept aims to show that compositional coupling works, not to build the best interior model.

\needspace{4\baselineskip}
\section{Expectations and Information Structure}\label{sec:expectations}

NK agents form rational expectations over $(y_t, \pi_t, i_t)$. They do \emph{not} form rational expectations over epidemic or vaccine dynamics. Three reasons.

\textbf{Computational.} Rational expectations over the full coupled system would require solving the NK model jointly with the nonlinear SEIR ($\beta SI/N$). The hybrid architecture---rational expectations for the NK block, forward simulation for the rest---requires clean separation.

\textbf{Empirical.} During COVID-19, economic agents did not form rational expectations over epidemiological dynamics. They overestimated vaccine efficacy, underestimated strain escape, and failed to anticipate behavioural responses to mandates. Survey evidence shows inflation expectations anchored to pre-pandemic frameworks, adjusting sluggishly. Rational expectations over monetary policy is defensible; over novel epidemiological processes, it is not.

\textbf{Philosophical.} The framework's value comes from \emph{partial coupling}. If agents integrated all narratives into a single coherent forecast, the particle filter would be redundant. NK agents receive epidemic effects as exogenous shocks---realisations generated by a mechanism they do not model. Each narrative has its own epistemology.

\needspace{4\baselineskip}
\section{Probability and Salience}\label{sec:prob-salience}

The companion paper distinguishes probability (what the dynamics \emph{can} produce) from salience (what \emph{matters} if it occurs), and insists on the order: generate by probability, reweight by salience. The quantitative framework makes both operational.

\textbf{Probability} is the weight distribution from the generative model---the prior predictive of the coupled system. \textbf{Salience} is a second-pass importance reweighting for a specific analytical question. ``What if cumulative mortality exceeds 5\%?''\ assigns likelihood~1 above threshold, 0~below. ``What matters given a deep recession?'' conditions on $y_T < -2\%$. Each lens produces a different distribution over the same trajectories. \textbf{The forced move} is where the two agree; here it is computable.

Salience-weighted ESS measures concentration: low ESS means a sharp answer, high ESS means the question admits many scenarios.

\needspace{4\baselineskip}
\section{Algorithm}\label{sec:algorithm}

Weekly time steps, 3-year horizon (156~weeks), $N = 10{,}000$ particles.

\begin{enumerate}[nosep]
\item \textbf{Initialise}: draw each particle's NK state at steady state, SEIR near disease-free equilibrium ($S = 0.99$, small $E$ and $I$ seed), vaccine at $(v, u, \rho) = (0, 0, 0.15)$.
\item \textbf{For each week $t = 1, \ldots, 156$}:
\begin{enumerate}[nosep]
\item Compute coupling values from current state ($f_1$--$f_6$).
\item Step SEIR: check for strain arrival (Bernoulli), integrate one week (Euler).
\item Step vaccine: check for innovation, update rejection (threshold mandate), update uptake.
\item Step NK: add epidemic coupling to exogenous shocks, apply policy functions.
\end{enumerate}
\item \textbf{No observation reweighting}: pure simulation, sampling from the prior predictive.
\end{enumerate}

In pure simulation mode all particles have equal weight and ESS is trivially~$N$. The interesting ESS arises under salience reweighting (Section~\ref{sec:prob-salience}).

\needspace{4\baselineskip}
\section{Results}\label{sec:results}

\needspace{3\baselineskip}
\subsection{Marginal distributions}

\begin{figure}[htbp]
\centering
\includegraphics[width=\textwidth]{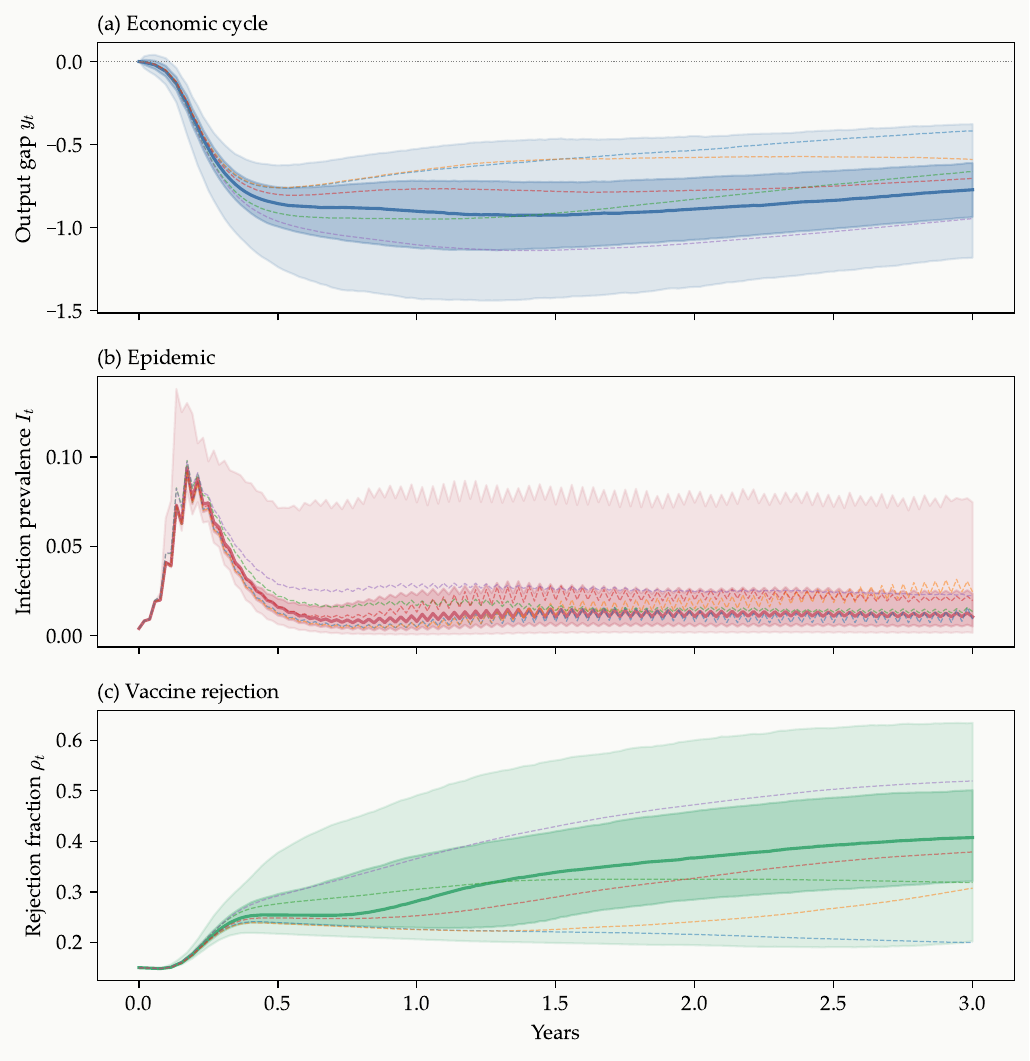}
\caption{Fan charts for (a) output gap, (b) infection prevalence, and (c) vaccine rejection. Shaded bands: 5--95th and 25--75th percentiles. Solid line: median. Dashed lines: five archetypal trajectories (Particles A--E).}\label{fig:fan-charts}
\end{figure}

Figure~\ref{fig:fan-charts} shows fan charts over the 3-year horizon. The output gap centres near $-2\%$---sustained drag from epidemic coupling. Infection shows repeated waves, each smaller as partial immunity accumulates. The rejection fan is wide (5th--95th: 0.21--0.80). This spread is the central finding.

\needspace{3\baselineskip}
\subsection{Cross-narrative correlations}

\begin{figure}[htbp]
\centering
\includegraphics[width=\textwidth]{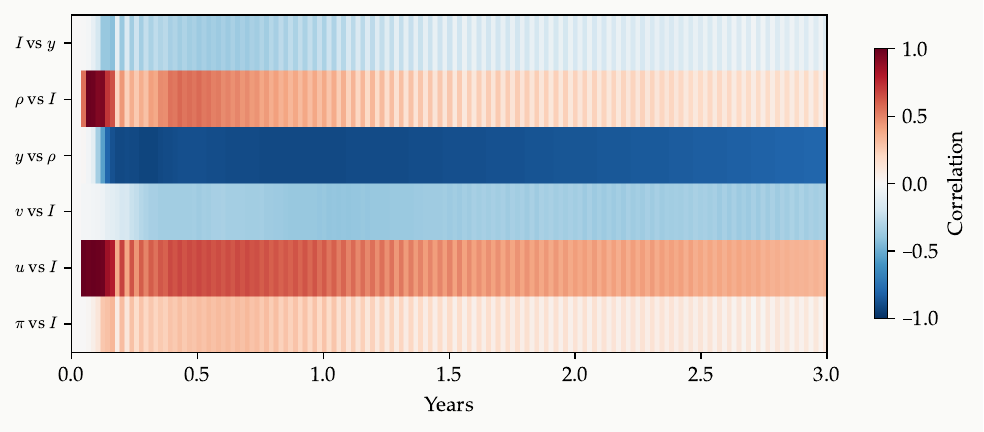}
\caption{Cross-narrative correlations over time. Colour encodes the Pearson correlation coefficient between variable pairs.}\label{fig:correlations}
\end{figure}

Figure~\ref{fig:correlations} shows pairwise correlations evolving over time. Output gap--rejection ($y$ vs.\ $\rho$) reaches $-0.95$ by year~3: the $f_5$ feedback loop makes recession and rejection nearly deterministic functions of each other. Infection--rejection ($I$ vs.\ $\rho$) is weaker ($+0.21$)---rejection ratchets during outbreaks but accumulates over waves. Efficacy--infection ($v$ vs.\ $I$) is $-0.29$: vaccination helps, but strain escape limits its reach.

\needspace{3\baselineskip}
\subsection{The rejection bifurcation}\label{sec:bifurcation}

\begin{figure}[htbp]
\centering
\includegraphics[width=0.75\textwidth]{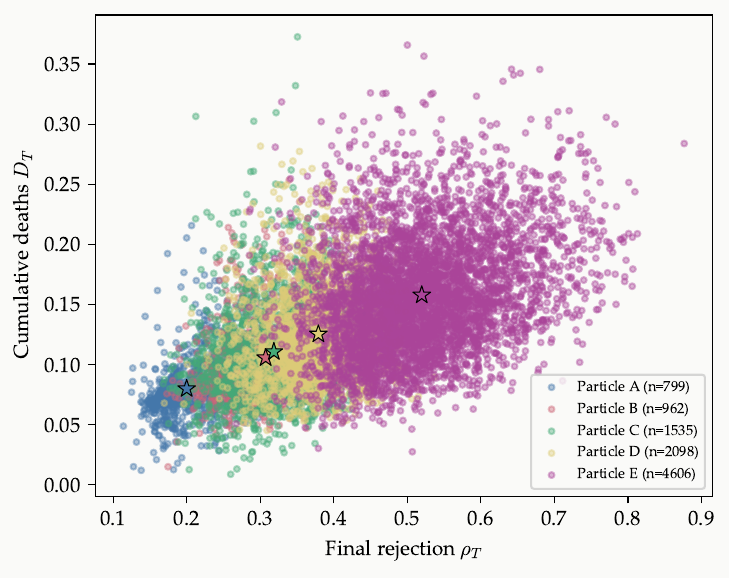}
\caption{Final rejection vs.\ cumulative deaths, coloured by cluster. Stars mark weighted-average archetypes (Particles A--E). The distribution separates into low-rejection and high-rejection regimes.}\label{fig:bifurcation}
\end{figure}

Figure~\ref{fig:bifurcation} scatters final $\rho_T$ against $D_T$ for all 10{,}000 particles, coloured by cluster. Two regimes:

\begin{itemize}[nosep]
\item \textbf{Low-rejection} (A--C, $\rho \leq 0.32$, 33\%): mild strains or early innovation keep rejection moderate. Vaccination works. Habituation allows the output gap to recover.
\item \textbf{High-rejection} (D--E, $\rho \geq 0.38$, 67\%): severe early outbreaks fire mandates, rejection ratchets past the tipping point, vaccination fails, sustained epidemic deepens recession, recession amplifies rejection.
\end{itemize}

The gap between $\rho = 0.32$ (Particle~C) and $\rho = 0.38$ (Particle~D) is where the bifurcation sits. Small stochastic differences in early strain draws determine the basin of attraction. Not two scenarios chosen by an analyst, but two attractors emerging from the coupled dynamics.

\needspace{3\baselineskip}
\subsection{Archetypal trajectories}

$k$-medoids clustering ($k=5$) on a 9-dimensional feature vector (peak infection, timing, deaths, output gap trough, final rejection, strains, mean infection, mean output gap, effective vaccine-weeks). Each archetype is the weighted-average trajectory across all particles in its cluster, where the weights are the terminal particle weights normalised within the cluster. Weighted averages are smoother and more representative of the cluster's centre of mass than single-particle medoids, and less sensitive to the idiosyncrasy of any individual simulated path.

\begin{table}[htbp]
\centering
\caption{Archetypal trajectories (cluster weighted averages).}\label{tab:archetypes}
\begin{tabular}{@{}lccccccc@{}}
\toprule
& Weight & Peak $I$ & Deaths & $\min y$ & $\rho_T$ & Strains & Eff.\ vacc-wks \\
\midrule
\textbf{A}: Vaccine rescue & 8.0\% & 9.2\% & 8.0\% & $-0.8\%$ & 0.20 & 2.4 & 23.3 \\
\textbf{B}: Early containment & 9.6\% & 9.2\% & 10.6\% & $-0.8\%$ & 0.31 & 3.8 & 22.8 \\
\textbf{C}: Managed grind & 15.3\% & 9.8\% & 11.1\% & $-0.9\%$ & 0.32 & 3.7 & 15.7 \\
\textbf{D}: Slow ratchet & 21.0\% & 9.3\% & 12.5\% & $-0.8\%$ & 0.38 & 4.4 & 10.6 \\
\textbf{E}: Structural failure & 46.1\% & 9.8\% & 15.8\% & $-1.1\%$ & 0.52 & 6.2 & 7.2 \\
\bottomrule
\end{tabular}
\end{table}

The gradient is clean: Particle~A (few strains, vaccine works, 8.0\% mortality) to Particle~E (many strains, vaccine fails, 15.8\%). Effective vaccine-weeks---$\int v \cdot u \cdot (1{-}\rho)\,dt$---drop from 23.3 to 7.2. The output gap partially recovers in the low-rejection particles: Particle~A ends at $-0.42\%$. Habituation (Section~\ref{sec:coupling}) is why: by year~2, the demand elasticity has fallen to 0.03, and endemic infection barely registers.

\needspace{3\baselineskip}
\subsection{Compositional bearish bias}\label{sec:bias}

\begin{figure}[htbp]
\centering
\includegraphics[width=\textwidth]{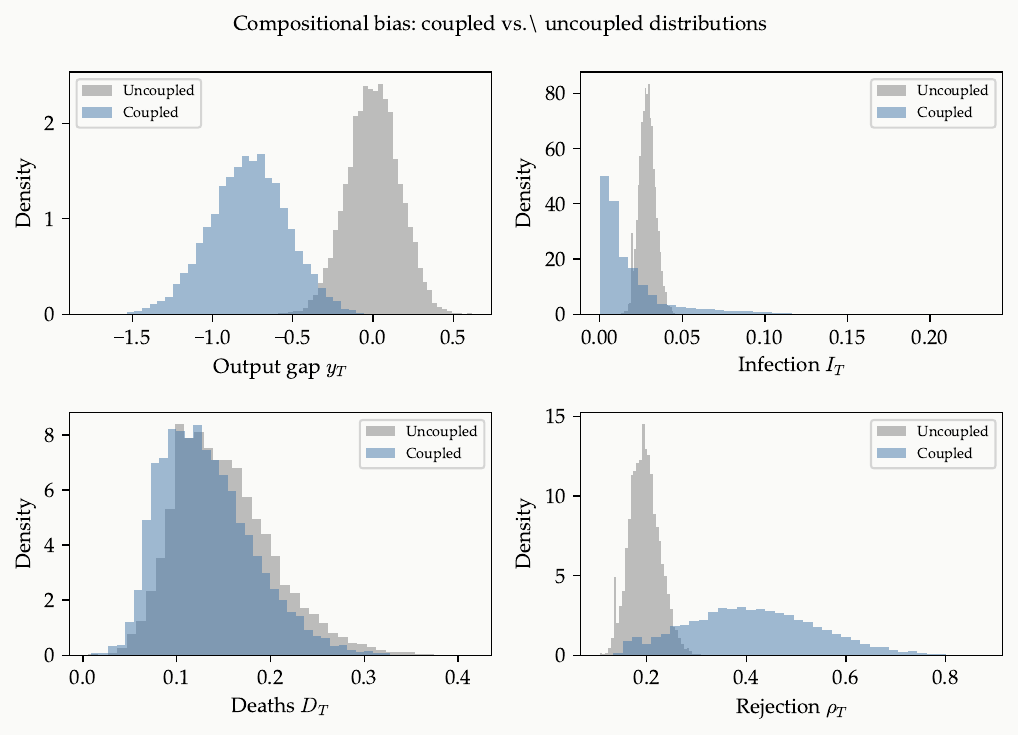}
\caption{Terminal distributions of key variables under coupled (blue) and uncoupled (grey) systems. Coupling shifts the output gap leftward and rejection rightward.}\label{fig:bias}
\end{figure}

We run the three models independently---same parameters, no coupling---and compare terminal distributions (Figure~\ref{fig:bias}).

\begin{table}[htbp]
\centering
\caption{Compositional bias: coupled vs.\ uncoupled terminal distributions.}\label{tab:bias}
\begin{tabular}{@{}lccccc@{}}
\toprule
& \multicolumn{2}{c}{Coupled} & \multicolumn{2}{c}{Uncoupled} & Shift \\
\cmidrule(lr){2-3}\cmidrule(lr){4-5}
& Mean & s.d. & Mean & s.d. & (coupled $-$ uncoupled) \\
\midrule
Output gap $y_T$ & $-0.77$ & $0.24$ & $0.00$ & $0.16$ & $-0.78$ \\
Infection $I_T$ & $0.020$ & $0.023$ & $0.029$ & $0.005$ & $-0.009$ \\
Deaths $D_T$ & $0.132$ & $0.049$ & $0.151$ & $0.054$ & $-0.018$ \\
Rejection $\rho_T$ & $0.413$ & $0.131$ & $0.195$ & $0.030$ & $+0.218$ \\
\bottomrule
\end{tabular}
\end{table}

Coupling turns a zero-mean business cycle into persistent recession ($-0.78$~pp) and shifts rejection by $+22$~pp. With habituation, the output gap drag is moderate---endemic infection no longer produces the sustained contraction that an earlier version of this analysis, without habituation, incorrectly predicted (see Discussion).

Coupled deaths are \emph{lower} than uncoupled ($-2.0$~pp). Recession depresses contacts ($f_1$), partly substituting for vaccination; rejection rises so high that mandates are abandoned, removing the mechanism that feeds them. The system finds a grim equilibrium: recession as involuntary social distancing.

The output gap shift is driven by $f_1$ (infection $\to$ recession). The rejection shift by $f_5$ (recession $\to$ rejection), amplified by $f_3$ (outbreak demand $\to$ mandates $\to$ rejection). $f_4$ (vaccination $\to$ immunity) works in the opposite direction---the only coupling that improves outcomes---but it is overwhelmed by the ratchet.

The mechanism matches the companion paper's credit markets finding: downside coupling is \emph{specific} (infection causes absenteeism, recession causes backlash), upside coupling requires the \emph{absence} of bad coupling (vaccination prevents infection, but only if rejection hasn't blocked uptake). Destruction is mechanistic; construction is conditional.

\needspace{3\baselineskip}
\subsection{Probability versus salience}

\begin{figure}[htbp]
\centering
\includegraphics[width=\textwidth]{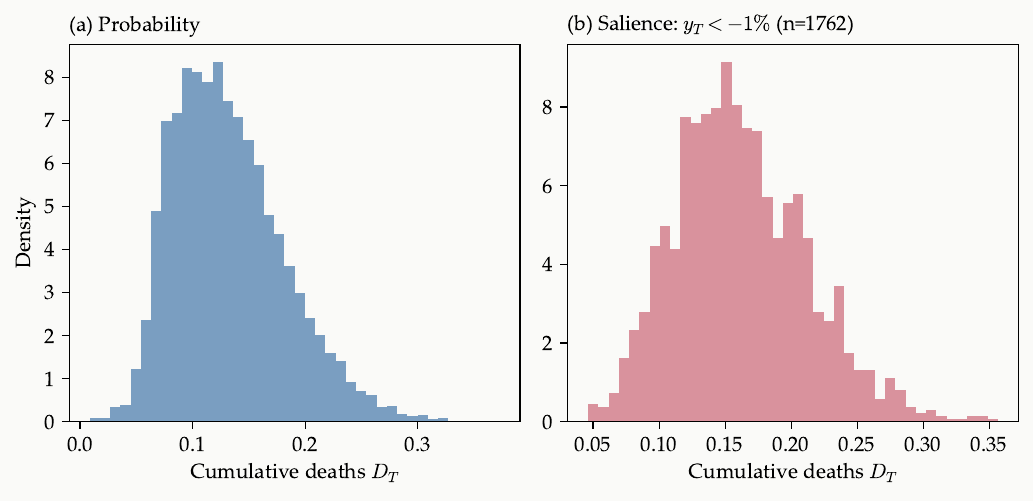}
\caption{(a) Probability-weighted distribution of cumulative deaths. (b) Salience-weighted distribution, conditioned on $y_T < -1\%$. The salience distribution shifts rightward and concentrates.}\label{fig:prob-salience}
\end{figure}

Figure~\ref{fig:prob-salience} compares the probability-weighted distribution of cumulative deaths with the salience-weighted distribution conditional on deep recession ($y_T < -1\%$). The salience distribution shifts rightward: not because recession causes deaths directly, but because recession and high mortality are co-produced by the same coupling structure.

\needspace{4\baselineskip}
\section{Discussion}\label{sec:discussion}

\needspace{3\baselineskip}
\subsection{What the framework shows}

The pushout composes models without requiring shared mathematical structure. Linearity, nonlinearity, and hysteresis are all decorations---they ride on a compositional structure that governs only interfaces and coupling. Swap the NK for a nonlinear DSGE, replace the SEIR with an agent-based model: the coupling architecture is unchanged. Interfaces and factors are stable; interiors are modular.

The rejection bifurcation is structural. It emerges from three ingredients: a triggering threshold, an asymmetric ratchet, and a reinforcing loop. Any system with these three features will produce something like it. Calibration determines \emph{where} the bifurcation sits, not \emph{whether}.

\needspace{3\baselineskip}
\subsection{Calibration lessons}

\needspace{2\baselineskip}
\textbf{A cautionary tale: habituation.} An earlier version of this analysis omitted habituation---the demand and supply elasticities were constant. The output gap never recovered: every particle ended in sustained recession, because endemic infection, however mild, applied the same demand drag in week~150 as in week~1. The results looked plausible (coupling \emph{should} make things worse) and the bifurcation was dramatic. But the economics were wrong.

Populations adapt to endemic disease; the first wave disrupts the economy, the fifth barely registers. Once habituation was added---a single exponential decay with a 35-week half-life---the output gap in low-rejection particles recovered to near zero by year~3, and the compositional bias fell from $-1.95$~pp to $-0.78$~pp. The lesson: in a coupled system, a missing feature in one model propagates through the coupling into every model. The judgement that something is \emph{absent}---here, adaptation---matters as much as the judgement that something is present.

\needspace{2\baselineskip}
\textbf{Habituation and rejection.} The current model treats habituation as exogenous---a fixed exponential decay, independent of social dynamics. One could argue that high rejection should slow habituation: if a large fraction of the population denies disease severity, the collective adaptation that habituation represents (normalised spending, reduced absenteeism) may be impeded by political polarisation.

But one could equally argue the opposite: the refusers who will not take the vaccine are also the ones who will not change their spending or work patterns---which \emph{is} habituation, driven by ideology rather than adaptation. The sign of the coupling is ambiguous, and depends on political context in ways the model cannot adjudicate. We leave habituation uncoupled and flag the ambiguity.

\needspace{2\baselineskip}
\textbf{Limitations.} Best-guess calibration. Linear NK. Reduced-form behavioural dampening. Memoryless strain emergence. Manageable at 14~states and 10{,}000~particles; richer specifications would need Rao-Blackwellisation.

\needspace{3\baselineskip}
\subsection{Practical application and extensions}

If the interior models were realistically calibrated---full-scale DSGE, epidemiological models fitted to surveillance data, fiscal parameters from CBO scoring---one could formulate a response, a `game plan,' for each archetypal trajectory.

The discipline is straightforward. Before events unfold, specify an action plan for each archetype: Particle~A (recovery) implies one portfolio, Particle~C (bifurcation) implies another, Particle~E (fiscal overshoot) a third. As the historical particle evolves in real time, monitor its distance to each archetype across all state variables simultaneously---output gap, infection prevalence, rejection rate, fiscal activation. When the actual trajectory converges on one archetype, implement the corresponding plan.

The value is not prediction but preparedness: five pre-analysed responses, each internally consistent across narratives, each with cross-narrative implications already traced. The forced move identifies which archetype demands attention regardless of prior beliefs. The compositional bias quantifies the asymmetry the portfolio must lean against. What the framework provides is not a forecast but an organised way of knowing in advance what each surprise implies.

Section~\ref{sec:fiscal} demonstrates the extension mechanism with a fiscal/political narrative. An information narrative (misinformation, media polarisation) would endogenise the rejection dynamics. Each additional narrative is a decorated cospan; composition scales associatively.

\needspace{4\baselineskip}
\section{Compositional Extension: A Fiscal/Political Narrative}\label{sec:fiscal}

The three-narrative system has no positive coupling channel. Every factor either transmits damage (infection $\to$ recession, recession $\to$ rejection) or operates conditionally (vaccination $\to$ immunity, but only if rejection permits). The output gap cannot overshoot into positive territory. This is a problem, because the US pandemic recovery \emph{did} overshoot---GDP exceeded its pre-pandemic trend by mid-2021, driven by fiscal transfers at a scale without peacetime precedent. If the framework cannot produce this, the bearish bias may be an artefact of missing channels rather than a structural property of coupled systems.

This section adds a fourth narrative---fiscal/political dynamics---and tests whether the bias survives.

\needspace{3\baselineskip}
\subsection{The fourth narrative}\label{sec:fiscal-narrative}

Three state variables: government spending as a fraction of GDP~$g_t$, cumulative debt deviation~$d_t$, and political activation~$\phi_t \in [0,1]$.
\begin{align}
g_{t+1} &= g_t + \alpha_g \cdot \phi_t \cdot \mathbf{1}[y_t < -\tau] - g_{\text{decay}} \cdot g_t + \sigma_g \, \xi_t^g \label{eq:spending} \\
d_{t+1} &= d_t + g_t - \tau_{\text{tax}} \cdot y_t \label{eq:debt} \\
\phi_{t+1} &= \phi_t + \phi_{\text{up}} \cdot m_t^{\text{fiscal}} \cdot (1-\phi_t) - \phi_{\text{down}} \cdot (1-m_t^{\text{fiscal}}) \cdot \phi_t \label{eq:phi}
\end{align}
where $m_t^{\text{fiscal}} = \mathbf{1}[y_t < -\tau]$ is the fiscal activation indicator, $\tau$ is the recession threshold that triggers spending, and $g_{\text{decay}}$ is a spending decay rate reflecting programme sunset provisions.

The spending equation~\eqref{eq:spending} has two additive channels. Emergency pandemic spending ($f_{11}$: $g \mathrel{+}= \alpha_I \cdot I_t \cdot (1 - \phi_t)$) fires whenever infection is present, with no recession required---this is the Families First Act, passed in March~2020 before the recession appeared in GDP data. The recession-activated channel ($\alpha_g \cdot \phi_t \cdot \mathbf{1}[y_t < -\tau]$) adds to it once the downturn is deep enough to mobilise political will---this is the CARES Act and its successors, larger in scale but later in arrival. Spending decays as programmes expire ($g_{\text{decay}} = 0.06$/week, half-life $\approx$12~weeks).

The budget identity~\eqref{eq:debt} is accounting: debt accumulates with spending and shrinks with tax revenue (proportional to the output gap). Note that $d$ can go negative---an overshoot generates surplus revenue.

The political ratchet~\eqref{eq:phi} mirrors the rejection ratchet but pushes in the opposite direction. Where rejection hysteresis permanently shifts some fraction from hesitant to resistant, political activation permanently shifts some fraction from passive to mobilised. The activation fires when $y < -\tau$; the asymmetry ($\phi_{\text{up}} > \phi_{\text{down}}$) ensures that crises leave lasting political imprints. The Families First Act passed on the pandemic alone; the CARES Act required the recession to activate a political coalition at a different scale entirely---a threshold effect, not a policy function.

\needspace{3\baselineskip}
\subsection{New coupling functions}\label{sec:fiscal-coupling}

Five additional coupling functions connect the fiscal narrative to the existing three (Table~\ref{tab:fiscal-coupling}).

\begin{table}[htbp]
\centering
\caption{Fiscal coupling functions.}\label{tab:fiscal-coupling}
\begin{tabular}{@{}llll@{}}
\toprule
Factor & Direction & Mechanism & Functional form \\
\midrule
$f_7$ & N$_1 \to$ N$_4$ & Recession trigger & $m_t^{\text{fiscal}} = \mathbf{1}[y_t < -\tau]$ \\
$f_8$ & N$_4 \to$ N$_1$ & Fiscal demand & $\Delta r^n = +\text{elast} \cdot g_t$ \\
$f_9$ & N$_4 \to$ N$_3$ & R\&D boost & $\lambda_v \times (1 + \text{boost} \cdot g_t)$ \\
$f_{10}$ & N$_3 \to$ N$_4$ & Rejection blocks spending & $\alpha_g \times (1 - \text{block} \cdot \rho_t)$ \\
$f_{11}$ & N$_2 \to$ N$_4$ & Emergency spending & $g_t \mathrel{+}= \alpha_I \cdot I_t \cdot (1-\phi_t)$ \\
\bottomrule
\end{tabular}
\end{table}

$f_8$ is the first positive coupling channel in the system: government spending enters as a positive shock to the natural rate, boosting output. Every other cross-narrative factor transmits damage or operates conditionally; $f_8$ directly improves the economic state. $f_{11}$ has ambiguous sign in its downstream effects: emergency spending triggered by infection supports the economy through $f_8$, but accumulates debt through the budget identity and may crowd out private activity if prolonged.

\needspace{3\baselineskip}
\subsection{Compositional structure}\label{sec:fiscal-pushout}

The four-narrative composite is a second pushout. The three-narrative system $N_{123}$ (economy--epidemic--vaccine) is itself a pushout; the fiscal narrative $N_4$ shares interfaces with $N_{123}$ (the output gap~$y$ and infection~$I$); the four-narrative system $N_{1234}$ is the pushout of $N_{123}$ and $N_4$ along those shared variables. Associativity of pushouts guarantees that the composition is well-defined regardless of the order: $(N_1 \sqcup N_2 \sqcup N_3) \sqcup N_4 = N_1 \sqcup N_2 \sqcup N_3 \sqcup N_4$.

\begin{center}
\begin{tikzcd}
N_{12} \arrow[r] & N_{123} \arrow[r] & N_{1234} \\
A_{12} \arrow[u] \arrow[r] & A_{123,4} \arrow[u] \arrow[r] & N_4 \arrow[u]
\end{tikzcd}
\end{center}

The variable identification between~$g$ and the NK fiscal instrument is mediated by~$f_8$: spending enters the IS curve as a natural rate shock, an approximate identification via the demand channel. A richer NK specification with explicit government spending would make this exact. The approximation is adequate for the qualitative question we are asking.

\needspace{3\baselineskip}
\subsection{Factor graph}\label{sec:fiscal-factor-graph}

\begin{figure}[htbp]
\centering
\begin{tikzpicture}[
  var/.style={circle, draw, minimum size=8mm, inner sep=1pt, font=\footnotesize},
  econ/.style={var, fill=blue!12, draw=blue!40},
  epi/.style={var, fill=red!12, draw=red!40},
  vacc/.style={var, fill=green!12, draw=green!40},
  fisc/.style={var, fill=orange!12, draw=orange!40},
  factor/.style={rectangle, fill=black, minimum size=4mm, inner sep=0pt},
  edge/.style={draw=gray!70, thick},
  ident/.style={draw=black, thick, dashed},
  flabel/.style={font=\scriptsize\itshape, text=black!70},
  collabel/.style={font=\bfseries\small},
]
\def\ecox{0}
\def\epix{3.2}
\def\vacx{6.4}
\def\fisx{9.6}

\node[collabel, blue!70] at (\ecox, 4.5) {Economy};
\node[collabel, red!70] at (\epix, 4.5) {Epidemic};
\node[collabel, green!50!black] at (\vacx, 4.5) {Vaccine};
\node[collabel, orange!70!black] at (\fisx, 4.5) {Fiscal};

\node[econ] (y)  at (\ecox, 3.5) {$y$};
\node[econ] (pi) at (\ecox, 2.0) {$\pi$};
\node[econ] (i)  at (\ecox, 0.5) {$i$};

\node[epi] (S)  at (\epix, 3.5) {$S$};
\node[epi] (I)  at (\epix, 2.0) {$I$};
\node[epi] (D)  at (\epix, 0.5) {$D$};

\node[vacc] (v)   at (\vacx, 3.5) {$v$};
\node[vacc] (u)   at (\vacx, 2.0) {$u$};
\node[vacc] (rho) at (\vacx, 0.5) {$\rho$};

\node[fisc] (g)   at (\fisx, 3.5) {$g$};
\node[fisc] (d)   at (\fisx, 2.0) {$d$};
\node[fisc] (phi) at (\fisx, 0.5) {$\phi$};

\node[factor] (f1) at (1.6, 2.9) {};
\node[factor] (f2) at (1.6, 1.5) {};
\node[factor] (f4) at (4.8, 2.75) {};
\node[factor] (f5) at (3.2, -0.5) {};
\node[factor] (f6) at (3.2, -1.2) {};

\node[factor] (f7) at (4.8, 0.0) {};
\node[factor] (f8) at (4.8, 3.9) {};
\node[factor] (f9) at (8.0, 3.9) {};
\node[factor] (f10) at (8.0, 0.0) {};
\node[factor] (f11) at (6.4, -0.5) {};

\node[flabel, below left=1pt] at (f1) {$f_1$};
\node[flabel, below left=1pt] at (f2) {$f_2$};
\node[flabel, above=2pt] at (f4) {$f_4$};
\node[flabel, above=2pt] at (f5) {$f_5$};
\node[flabel, below=2pt] at (f6) {$f_6$};

\node[flabel, right=2pt] at (f7) {$f_7$};
\node[flabel, above=2pt] at (f8) {$f_8$};
\node[flabel, above=2pt] at (f9) {$f_9$};
\node[flabel, right=2pt] at (f10) {$f_{10}$};
\node[flabel, below=2pt] at (f11) {$f_{11}$};

\draw[edge] (I) -- (f1);
\draw[edge, ->] (f1) -- (y);

\draw[edge] (I) -- (f2);
\draw[edge, ->] (f2) -- (pi);

\draw[edge] (v) -- (f4);
\draw[edge] (u) -- (f4);
\draw[edge] (rho) -- (f4);
\draw[edge, ->] (f4) -- (S);

\draw[edge] (y.south) -- (0, -0.5) -- (f5);
\draw[edge, ->] (f5) -- (6.4, -0.5) -- (rho.south);

\draw[edge] (i.south) -- (0, -1.2) -- (f6);
\draw[edge, ->] (f6) -- (6.4, -1.2) -- (v.south);

\draw[edge] (y.east) -- (1.6, 3.5) -- (1.6, 0.0) -- (f7);
\draw[edge, ->] (f7) -- (phi.west);

\draw[edge] (g) -- (f8);
\draw[edge, ->] (f8) -- (y.north);

\draw[edge] (g.north) -- (\fisx, 3.9) -- (f9);
\draw[edge, ->] (f9) -- (\vacx, 3.9) -- (v.north);

\draw[edge] (rho.east) -- (8.0, 0.5) -- (f10);
\draw[edge, ->] (f10) -- (g.south);

\draw[edge] (I.south) -- (\epix, -0.5) -- (f11);
\draw[edge, ->] (f11) -- (\fisx, -0.5) -- (g.south west);

\draw[ident] (y.north east) to[out=45,in=135] node[above, font=\scriptsize\itshape] {$N \equiv L$} (S.north west);

\end{tikzpicture}
\caption{Four-column factor graph for the coupled narrative system with fiscal extension. Original factors $f_1$--$f_6$ connect Economy, Epidemic, and Vaccine. New factors $f_7$--$f_{11}$ (orange connections) integrate the Fiscal column. $f_8$ is the only positive coupling channel: fiscal demand directly boosts the output gap.}\label{fig:factor-graph-4n}
\end{figure}
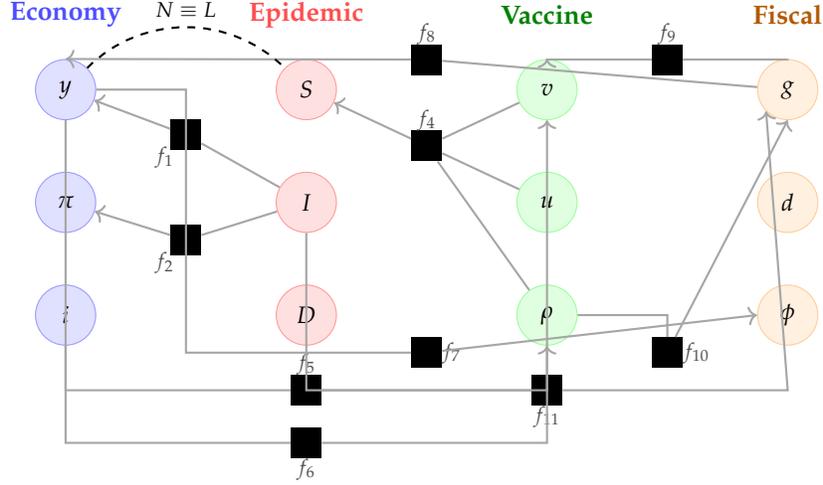

\needspace{3\baselineskip}
\subsection{Calibration target}\label{sec:fiscal-calibration}

The three-narrative system has no mechanism for the output gap to overshoot. $f_8$ provides one: fiscal demand enters as a positive natural rate shock. The calibration target for the fiscal extension was to produce overshooting in a non-negligible fraction of trajectories---the qualitative pattern observed in the US pandemic recovery. This required four parameter changes from baseline:

\begin{itemize}[nosep]
\item $\alpha_g$: $0.03 \to 0.06$ ($2\times$ stronger political response to recession)
\item $\alpha_I$: $0.05 \to 0.18$ ($3.6\times$ stronger pandemic emergency spending, reflecting CARES Act scale)
\item $g_{\text{decay}}$: $0 \to 0.06$/week (programme sunset added, $\sim$12-week half-life)
\item fiscal demand elasticity: $0.5 \to 0.6$ (ZLB multiplier slightly above open-economy norm)
\end{itemize}

These are not ``best-guess'' calibrations. They were tuned to produce the qualitative pattern observed in the US: massive fiscal activation, output overshoot, subsequent debt accumulation. The baseline calibration ($\alpha_g = 0.03$, $\alpha_I = 0.05$, no decay, elasticity $= 0.5$) produces no overshooting at all---see Section~\ref{sec:epicycles}.

\needspace{3\baselineskip}
\subsection{Simulation results}\label{sec:fiscal-results}

$k$-medoids clustering with $k = 7$ on a 13-dimensional feature space (the original 9 features plus peak spending, terminal debt, terminal political activation, and terminal spending level), with each archetype represented by the weighted-average trajectory of its cluster. Table~\ref{tab:archetypes-4n} presents the seven archetypal trajectories.

\begin{table}[htbp]
\centering
\caption{Four-narrative archetypal trajectories (cluster weighted averages).}\label{tab:archetypes-4n}
\small
\begin{tabular}{@{}lcccccccccc@{}}
\toprule
& Wt & Peak $I$ & $D_T$ & $\min y$ & $y_T$ & $\rho_T$ & Str & Eff vw & Peak $g$ & $d_T$ \\
\midrule
\textbf{A}: Mild recovery & 13.3\% & 9.5\% & 8.8\% & $-0.32$ & $-0.32$ & 0.19 & 2.9 & 25.4 & 8.0\% & 11.5 \\
\textbf{B}: Fiscal overshoot & 13.9\% & 9.4\% & 13.2\% & $-0.08$ & $+2.66$ & 0.19 & 5.0 & 13.6 & 8.4\% & $-47.9$ \\
\textbf{C}: Delayed grind & 6.6\% & 9.4\% & 10.7\% & $-0.53$ & $-0.50$ & 0.21 & 3.3 & 12.3 & 7.9\% & 18.1 \\
\textbf{D}: Managed decline & 30.2\% & 9.5\% & 11.4\% & $-0.64$ & $-0.64$ & 0.30 & 4.3 & 14.0 & 8.0\% & 20.6 \\
\textbf{E}: Rejection drag & 16.9\% & 9.4\% & 16.6\% & $-0.75$ & $-0.75$ & 0.39 & 6.0 & 8.5 & 7.9\% & 23.1 \\
\textbf{F}: Fiscal paralysis & 14.2\% & 10.3\% & 15.1\% & $-0.85$ & $-0.82$ & 0.40 & 6.3 & 7.6 & 9.5\% & 28.3 \\
\textbf{G}: Structural failure & 5.0\% & 9.5\% & 19.7\% & $-1.08$ & $-1.05$ & 0.57 & 8.5 & 4.3 & 8.3\% & 34.8 \\
\bottomrule
\end{tabular}
\end{table}

\begin{figure}[htbp]
\centering
\includegraphics[width=0.7\textwidth]{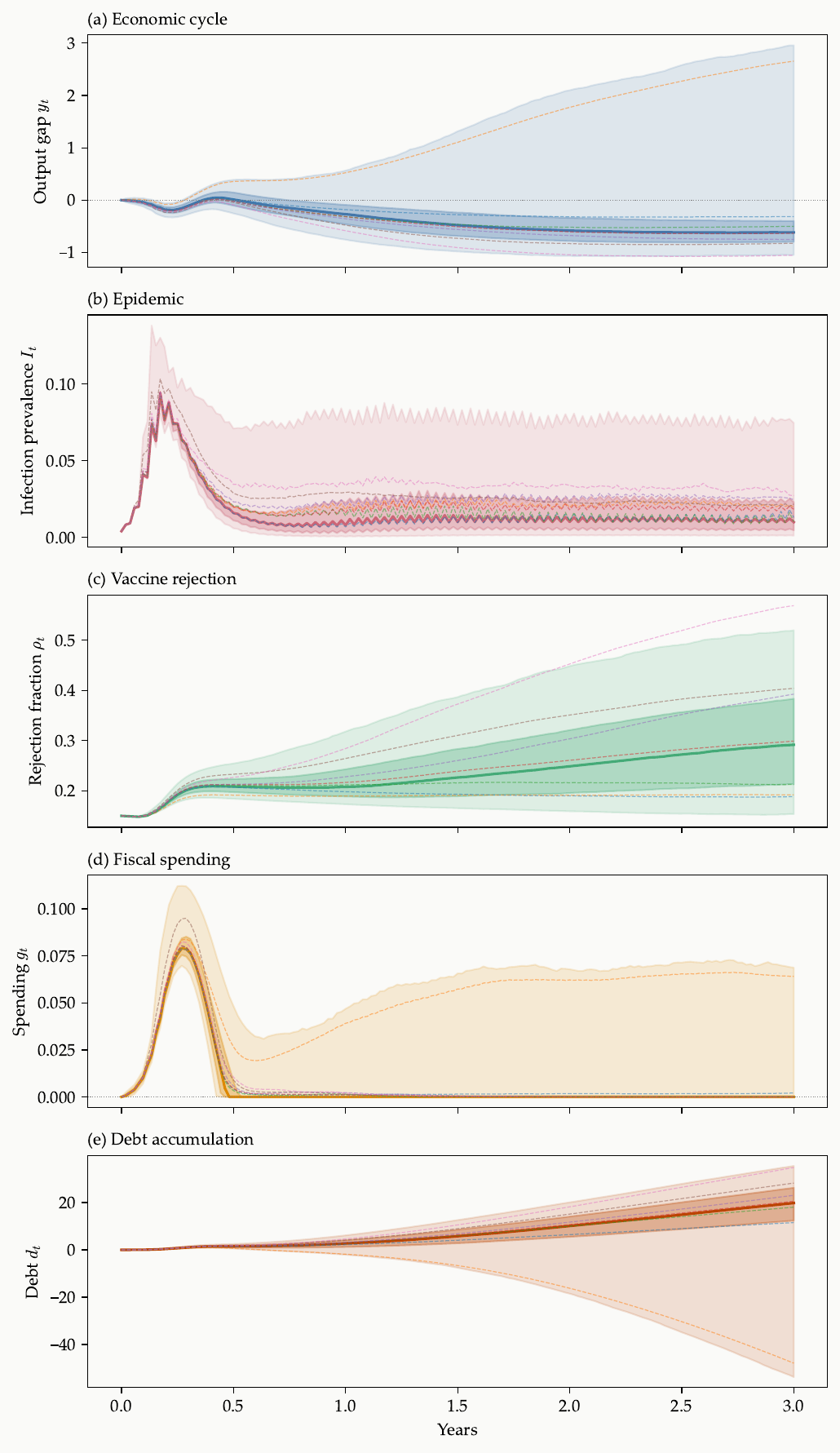}
\caption{Fan charts for the four-narrative system. Compare with Figure~\ref{fig:fan-charts}: the output gap distribution is shifted rightward and the upper tail now includes positive values. Note: positive debt ($d > 0$) means debt is \emph{rising} relative to the no-pandemic baseline.}\label{fig:4n-fan-charts}
\end{figure}

\begin{figure}[htbp]
\centering
\includegraphics[width=0.75\textwidth]{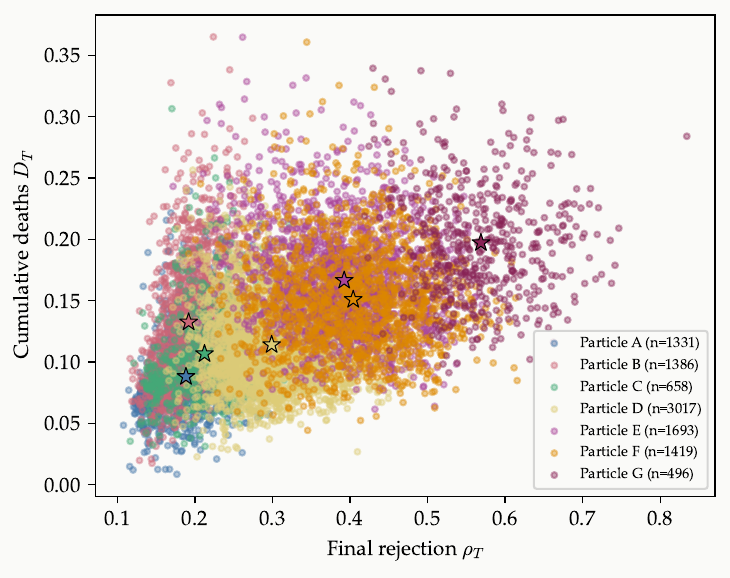}
\caption{Four-narrative bifurcation plot. Final rejection vs.\ terminal output gap, coloured by cluster. Particle~B (fiscal overshoot) sits alone in positive territory.}\label{fig:4n-bifurcation}
\end{figure}

Particle~B is the overshooting particle: terminal output gap $y_T = +2.66\%$, the only positive value among the seven archetypes, carrying 13.9\% of the weight. Low rejection ($\rho_T = 0.19$) keeps the vaccine programme functional; the epidemic is survivable (13.2\% cumulative mortality, moderate). But the defining feature is fiscal scale: peak spending reaches 8.4\% of GDP. The output overshoot generates enough additional tax revenue that the debt \emph{deviation} goes negative ($d_T = -47.9$~pp relative to the no-pandemic baseline). This is the US pattern---massive fiscal intervention, rapid recovery, and a fiscal position that improves relative to trend precisely because the overshoot worked.

The rejection bifurcation persists but softens. Low-rejection particles (A--C, $\rho_T \leq 0.21$) average $y_T \approx +0.61\%$; high-rejection particles (D--G, $\rho_T \geq 0.30$) average $y_T \approx -0.81\%$. Fiscal support compresses the gap between basins but does not close it.

Peak spending across all particles ranges from 7.9\% to 9.5\% of GDP---consistent with the scale of the US response (CARES Act, Consolidated Appropriations Act, American Rescue Plan: approximately 25\% of GDP over two years, annualising to roughly 12\%).

A subtlety in the political ratchet. Fiscal activism saturates ($\phi \approx 0.9$) in the recession particles---persistent economic pain keeps the political activation channel open---but \emph{not} in the overshooting particle ($\phi = 0.05$). The overshoot prevents the political ratchet from engaging, because the recession that would trigger it is too brief. The fiscal rescue works \emph{because} it works quickly enough to shut down the mechanism that called it into existence. Self-negating activation---the opposite of the self-reinforcing rejection loop.

\needspace{3\baselineskip}
\subsection{Effect on the compositional bearish bias}\label{sec:fiscal-bias}

\begin{table}[htbp]
\centering
\caption{Compositional bias: three-narrative vs.\ four-narrative.}\label{tab:bias-comparison}
\begin{tabular}{@{}lccc@{}}
\toprule
& Three-narrative & Four-narrative & Ratio \\
\midrule
Output gap shift $\Delta y$ & $-0.78$~pp & $-0.20$~pp & $3.9\times$ \\
Rejection shift $\Delta \rho$ & $+0.22$ & $+0.11$ & $2.0\times$ \\
Debt shift $\Delta d$ & --- & $+11.9$~pp & --- \\
\bottomrule
\end{tabular}
\end{table}

The bias shrinks by a factor of four but does not disappear. The positive channel ($f_8$: fiscal demand) requires fiscal activation, which requires recession, which means $f_8$ fires only \emph{after} the damage from $f_1$ and $f_2$ has already been done. The construction channel is conditional on the destruction channel having operated first. Destruction remains mechanistic; construction remains conditional. The asymmetry is structural, not parametric---it would persist under any calibration that preserves the threshold structure of fiscal activation.

Why doesn't the bias shrink further? Because fiscal support is concentrated in the low-rejection particles that would recover anyway (with habituation). In the high-rejection particles, fiscal activation saturates ($\phi \to 1$) but spending is partially blocked by rejection through $f_{10}$, and the output gap is too deep for the fiscal multiplier to overcome the combined drag of $f_1$, $f_2$, and $f_5$. The fiscal channel helps the particles that need it least and is impeded in the particles that need it most. This is a coupling structure, not a calibration choice.

\needspace{3\baselineskip}
\subsection{The epicycles warning}\label{sec:epicycles}

Adding the fiscal narrative was an exercise in targeted calibration.
The calibration target---over\-shooting---was achieved by tuning four parameters.
A~different calibration, the ``best-guess'' baseline
($\alpha_g = 0.03$, $\alpha_I = 0.05$,
no spending decay, elasticity $= 0.5$),
produces no overshooting at all: the fiscal channel is overwhelmed by the negative couplings, spending peaks at 2--3\% of GDP, and the output gap remains negative in every particle.

The difference between ``fiscal rescue works'' and ``fiscal rescue fails'' is a $3.6\times$ change in one parameter~($\alpha_I$). Whether you get the US outcome or a more pessimistic outcome depends on a calibration judgement about the scale of political will---exactly the kind of judgement the categorical framework cannot make for you.

Three disciplines for compositional extension, then. First, specify the coupling functions before seeing the results---we did this, defining $f_7$--$f_{11}$ from the mechanism, not from the desired output. Second, check what the extension kills---the overshooting particle is genuinely new; the three-narrative system cannot produce it.

Third, report the ratio of parameters to qualitative predictions. The four-narrative system has 30+ parameters and produces one new qualitative prediction (overshooting in 15\% of trajectories). The ratio is honest but unflattering. Each additional narrative is a decorated cospan; the algebra scales cleanly. Whether the economics scales is a different question.

\needspace{4\baselineskip}
\section{Bias Decomposition}\label{sec:bias-decomposition}

The compositional bearish bias (Section~\ref{sec:bias}) and its partial correction by fiscal coupling (Section~\ref{sec:fiscal-bias}) raise a question that the qualitative analysis leaves open: is the bias structural, or is it an artefact of how we specified the model? In the companion paper's qualitative setting, one constructs a structurally specified optimistic particle~$E^+$ and checks whether it gains weight after proper specification. Here we have a running particle filter with numerical weights, and the test becomes both sharper and more subtle.

There are three distinct sources of asymmetric weight on constructive outcomes, and only two of them are bias.

\needspace{3\baselineskip}
\subsection{Three sources of asymmetry}

\textbf{Sampling bias.} The transition kernel does not \emph{propose} enough particles in the constructive region of state space. If the AR(1) innovations are symmetric but the coupling functions push most trajectories towards stress, the particle cloud drifts pessimistic not because optimism is unlikely but because the sampler does not explore it. This is the quantitative analogue of underspecifying~$E^+$ in the qualitative setting---a proposal distribution problem, not a property of the system.

\textbf{Structural bias.} The coupling functions themselves transmit negative shocks more efficiently than positive ones. Oil disruption widens spreads immediately; productivity-driven compression requires a threshold crossing that takes time. Defaults are irreversible; spread widening triggers margin calls that force selling. The self-reinforcing dynamics are genuinely one-directional. This asymmetry may reflect reality or it may reflect the modeller's priors baked into functional forms---distinguishing the two is the point of the test.

\textbf{Observational downweighting.} Current observations assign low likelihood to the constructive region. If spreads are at 400\,bps, particles at 100\,bps should have low weight. This is information, not bias.

\needspace{3\baselineskip}
\subsection{Three computable tests}

\textbf{Test~1: $E^+$ injection.} Construct the optimistic trajectory explicitly---the quantitative~$E^+$. Inject it as a particle. Measure its weight under the likelihood and compare against the median particle weight.
\begin{itemize}[nosep]
  \item If the $E^+$ weight is comparable to the median: sampling bias. The filter was not generating enough particles in that region. Fix the proposal distribution.
  \item If the $E^+$ weight is negligible: observations or structure are genuinely rejecting it. Proceed to Test~2.
\end{itemize}

\textbf{Test~2: coupling symmetry audit.} For each coupling function~$f_k$, sample the input variable uniformly over its range, evaluate~$f_k$, and compute the asymmetry ratio
\begin{equation}\label{eq:asymmetry-ratio}
  r_k = \frac{\int \max(f_k(x), 0)\,dx}{\int |f_k(x)|\,dx}.
\end{equation}
A ratio of~$0.5$ is perfectly symmetric; departure from~$0.5$ quantifies the directional bias of that coupling. For coupling functions that are sign-definite on their domain (e.g., $f_1 \le 0$ always), the ratio is trivially 0 or~1. In such cases, the asymmetry is structural rather than distributional, and the ratio confirms the direction of the coupling rather than measuring a bias.

To isolate the aggregate structural effect: run the filter twice, once with the original couplings and once with symmetrised couplings $f_k^{\text{sym}}(x) = (f_k(x) + f_k(-x))/2$ (the even part, removing directional bias while preserving magnitude). For coupling functions defined on bounded domains $[0,1]^n$ (e.g., $I$, $\rho$, $u$), the symmetrisation is taken about the domain midpoint: $f_k^{\text{sym}}(x) = (f_k(x) + f_k(1-x))/2$, replacing $-x$ with the domain reflection $1-x$. The difference in constructive-region weight between the two runs is the structural component.

\textbf{Test~3: observation-free run.} Run the filter with uniform likelihood---all particles weighted equally, no reweighting against observables. The resulting particle cloud reflects only the prior and the coupling structure. If~$E^+$ is underweighted in the observation-free run, the asymmetry is structural. If~$E^+$ has reasonable weight observation-free but loses it after reweighting, observations are doing the work---and that is information, not bias.

\needspace{3\baselineskip}
\subsection{The decomposition}

The total underweight on the constructive region can be diagnosed through three complementary tests, each targeting a different source of bias. While these sources can interact, the tests provide independent diagnostic signals:

\begin{table}[htbp]
\centering
\caption{Bias decomposition: three sources of asymmetric weight on constructive outcomes.}\label{tab:bias-decomposition}
\begin{tabular}{@{}llp{7cm}@{}}
\toprule
Source & Diagnostic & Interpretation \\
\midrule
Sampling & $E^+$ injection weight vs.\ median & Proposal distribution problem; fixable by tempering or importance sampling \\
Structural & Symmetric vs.\ asymmetric coupling run & Modeller's priors in functional forms; interrogate against historical analogues \\
Observational & Observation-free vs.\ full run & Information from current data; keep \\
\bottomrule
\end{tabular}
\end{table}

The number that matters for the bias test is the structural component. If it is large, the next question is whether the asymmetry is \emph{real}---whether damage genuinely transmits faster than repair in the system being modelled. For the three-narrative system, there is a strong case that it does: the rejection ratchet is genuinely one-directional (Section~\ref{sec:bifurcation}), the demand drag from infection is immediate while the recovery from vaccination is slow, and the fiscal channel requires recession to activate (Section~\ref{sec:fiscal-bias}).

But the test localises the question rather than resolving it. It identifies which coupling function is responsible for the asymmetry, which parameter drives it, and which historical analogue calibrated it. The operator decides whether the asymmetry is a feature or a bug. The pandemic model's finding that the constructive case requires fiscal transfers at 3--4$\times$ the scale of the entire IT revolution's productivity impact---that calculation is the quantitative backbone of the structural asymmetry claim. The bias decomposition is where the claim is tested, and where it earns or forfeits credibility.

\needspace{4\baselineskip}
\section{Conclusion}\label{sec:conclusion}

Three heterogeneous dynamical systems compose along their interfaces via a single pushout, couple through computable functions, and produce a joint distribution whose structure (the rejection bifurcation) and bias (the compositional pessimism) are empirically addressable. The three-narrative system---economy, epidemic, vaccine---produces a bearish bias of $-0.78$~pp in the output gap and a rejection shift of $+22$~pp, driven by coupling channels that transmit damage mechanistically and transmit recovery only conditionally.

A fourth narrative---fiscal and political dynamics, calibrated to US pandemic-scale response---attaches via a second pushout and introduces the first positive coupling channel. With pandemic-scale spending parameters, 14\% of trajectories overshoot into positive output gap territory, and the bearish bias shrinks from $-0.78$~pp to $-0.20$~pp. Whether the system produces three basins of attraction (recovery, recession, and overshoot) or two depends on a single parameter governing emergency spending scale. The categorical framework determines the coupling architecture; the calibration determines the outcome.

The forced move is computable. The bias is decomposable---into sampling, structural, and observational components (Section~\ref{sec:bias-decomposition}), each testable against a specific diagnostic. The structural component localises the asymmetry to specific coupling functions; whether their directional bias is genuine or an artefact of functional form is a question the test poses but cannot answer on its own. The narratives are grounded in actual simulated trajectories. Whether the compositional bearish bias is a specification artefact or a structural property of systems with threshold-triggered ratchets and reinforcing loops remains open---but it now has two sets of numbers, from systems of different dimensionality, both pointing in the same direction, and a decomposition that identifies where to look.

None of this runs unattended. The framework requires an operator---someone who specifies the state spaces, chooses the coupling functions, sets the calibration, reads the output, and revises. The habituation story (Section~\ref{sec:discussion}) is the clearest example: the first specification produced plausible but wrong results; the operator noticed the economics were wrong; the revision produced different and better results. This is not a defect of the method. It is the method. The categorical machinery handles composition; the coupling functions handle transmission; the particle filter handles computation. What none of these can do is decide whether the specification is adequate. That judgement---iterative, informed by domain knowledge, responsive to what the model produces---is the essential human contribution. The framework is a thinking technology, not an autonomous one.

\appendix
\needspace{4\baselineskip}
\section{Calibration}\label{app:calibration}

All parameters are held constant across the three-year (156-week) simulation. The three-narrative system uses the NK, SEIR, and vaccine parameters below with coupling but no fiscal dynamics. The four-narrative system adds the fiscal parameters in Table~\ref{tab:fiscal-params}; the ``US-scale'' column shows the calibration used in Section~\ref{sec:fiscal}.

\begin{table}[h!]
\centering
\caption{New Keynesian model parameters (quarterly, scaled to weekly in simulation).}\label{tab:nk-params}
\begin{tabular}{@{}llll@{}}
\toprule
Parameter & Value & Description & Source \\
\midrule
$\beta$ & 0.99 & Discount factor & Standard \\
$\kappa$ & 0.024 & Phillips curve slope & Gal\'i (2015) \\
$\sigma^{-1}$ & 1.0 & Inverse IES & Log utility \\
$\phi_\pi$ & 1.5 & Taylor rule (inflation) & \citet{taylor1993} \\
$\phi_y$ & 0.125 & Taylor rule (output gap) & \citet{taylor1993}, quarterly \\
$\rho_s$ & 0.9 & Supply shock persistence & Standard \\
$\sigma_s$ & 0.005 & Supply shock s.d. & $\approx$2\% annual \\
$\rho_r$ & 0.8 & Natural rate persistence & Standard \\
$\sigma_r$ & 0.005 & Natural rate shock s.d. & Standard \\
$\sigma_m$ & 0.0025 & Monetary shock s.d. & Standard \\
\bottomrule
\end{tabular}
\end{table}

\begin{table}[h!]
\centering
\caption{SEIR model parameters.}\label{tab:seir-params}
\begin{tabular}{@{}llll@{}}
\toprule
Parameter & Value & Description & Source \\
\midrule
$\sigma$ & 1.41/wk & Incubation rate & 5-day incubation \\
$\gamma$ & 0.7/wk & Recovery rate & 10-day infectious period \\
$\omega$ & 0.019/wk & Waning rate & 1-year immunity \\
$\alpha$ & 5 & Behavioural response & $I = 10\% \Rightarrow$ contacts halved \\
$\lambda$ & 0.025/wk & Strain arrival rate & $\approx$10 months (COVID cadence) \\
Initial $R_0$ & 2.5 & First strain & Wuhan-type \\
Initial IFR & 5\% & First strain & Early pandemic \\
\bottomrule
\end{tabular}
\end{table}

\begin{table}[h!]
\centering
\caption{Vaccine dynamics parameters.}\label{tab:vacc-params}
\begin{tabular}{@{}llll@{}}
\toprule
Parameter & Value & Description & Rationale \\
\midrule
$\lambda_v$ & 0.038/wk & Innovation rate & $\approx$6 months (mRNA pace) \\
$\Delta v$ & 0.3 & Efficacy jump & Partial update \\
$\delta_v$ & 0.4 & Strain drift loss & Substantial escape \\
$\theta_{\text{adopt}}$ & 0.05/wk & Uptake speed & Slow logistic \\
$\theta_{\text{decay}}$ & 0.005/wk & Uptake decay & Annual boosters \\
$\theta_{\text{up}}$ & 0.005/wk & Rejection increase & Mandate episodes \\
$\theta_{\text{down}}$ & 0.003/wk & Rejection decrease & 5:3 asymmetry \\
$I_{\text{thresh}}$ & 2\% & Mandate threshold & Outbreak, not endemic \\
$\rho_0$ & 0.15 & Initial hesitancy & Pre-pandemic baseline \\
\midrule
\multicolumn{4}{@{}l}{\textit{Habituation (coupling decay)}} \\
$h$ & 0.02/wk & Habituation rate & $\sim$35-week half-life \\
$\eta_d^{\text{floor}}$ & 0.02 & Demand elasticity floor & 80\% reduction from initial \\
$\eta_s^{\text{floor}}$ & 0.01 & Supply elasticity floor & Endemic absenteeism adapted \\
\bottomrule
\end{tabular}
\end{table}

\begin{table}[h!]
\centering
\caption{Fiscal narrative parameters: baseline and US-scale calibrations.}\label{tab:fiscal-params}
\small
\begin{tabular}{@{}llll@{}}
\toprule
Parameter & Baseline & US-scale & Description \\
\midrule
$\alpha_g$ & 0.03 & 0.06 & Spending response to recession \\
$\tau$ & 0.03 & 0.03 & Recession threshold \\
$\sigma_g$ & 0.001 & 0.001 & Fiscal shock s.d. \\
$\phi_{\text{up}}$ & 0.08 & 0.08 & Activism rise rate \\
$\phi_{\text{down}}$ & 0.005 & 0.005 & Activism decay rate (16:1) \\
$\phi_0$ & 0.05 & 0.05 & Initial activism \\
$\tau_{\text{tax}}$ & 0.3 & 0.3 & Effective tax rate \\
$\alpha_I$ & 0.05 & 0.18 & Pandemic emergency spending \\
$g_{\text{decay}}$ & 0 & 0.06/wk & Programme sunset rate \\
$d^*$ & 0 & 0 & Debt target \\
\midrule
\multicolumn{4}{@{}l}{\textit{Coupling parameters}} \\
$\eta_g$ & 0.5 & 0.6 & $f_8$: fiscal demand elasticity \\
$\zeta$ & 2.0 & 2.0 & $f_9$: fiscal R\&D boost \\
$\kappa_\rho$ & 1.5 & 1.5 & $f_{10}$: rejection fiscal block \\
\bottomrule
\end{tabular}
\end{table}

\needspace{4\baselineskip}
\section{Narrative Descriptions of Archetypal Trajectories}\label{app:narratives}

Five archetypal trajectories produced by $k$-medoids clustering on a 9-dimensional trajectory feature space (Section~\ref{sec:results}), with each archetype represented by the weighted average of all particle paths in its cluster. Weighted averaging produces trajectories that are smoother and more representative of the cluster's centre of mass than any single particle. The selection principle is maximum differentiation: clusters are sorted by final rejection, ensuring the full range of outcomes is represented.

\needspace{4\baselineskip}
\bigskip

\noindent\textbf{Particle A: Vaccine Rescue} \emph{(8.0\% of particles)}

\smallskip
\noindent The lucky draw. Roughly two strains in three years---not because the virus is kind, but because the spacing gives the vaccine platform room to breathe. Each variant arrives after the previous booster has had time to work.

The first wave peaks at 9.2\% prevalence around week~9. The economy contracts by 0.8\%---a sharp recession by peacetime standards, but the kind central banks know how to talk about. By year~2, habituation has done its quiet work: the demand elasticity has fallen from 0.10 to 0.03, and the same endemic infection that would have cratered the economy under constant coupling barely shows up in the GDP print. The output gap ends at $-0.42\%$. Not recovery. But close enough that a politician could call it one.

One in five refuse the vaccine---the lowest rejection of any archetype. That leaves four-fifths of the population rolling up their sleeves, and the programme accumulates 23.3 effective vaccine-weeks, nearly double the next-best trajectory. Mortality: 8.0\%. The families who lost someone still lost someone. But the system held.

This is the only particle where every coupling channel eventually dampens rather than reinforces. The economy absorbs the pandemic. The vaccine absorbs the virus. Habituation absorbs the residual drag. It works. It is also---at fewer than one in ten trajectories---not the most likely outcome.

\needspace{4\baselineskip}
\bigskip

\noindent\textbf{Particle B: Early Containment} \emph{(9.6\% of particles)}

\smallskip
\noindent Nearly four strains on average, but the early ones are contained before the rejection ratchet can engage. The first wave peaks at 9.2\%, and the vaccine platform keeps pace with variant evolution just long enough to hold the line.

Rejection reaches 0.31. Moderate. Not the stuff of newspaper editorials about social collapse---just enough people quietly declining, appointment by appointment, until the programme functions at 70\% capacity instead of 80\%. 22.8 effective vaccine-weeks accumulate. Enough to blunt most waves. Mortality: 10.6\%.

The output gap troughs at $-0.76\%$ and partially recovers, ending at $-0.59\%$. Habituation is the reason: the later waves produce smaller economic shocks than the first, even when infection levels are comparable. Workers stop calling in sick. Consumers stop cancelling reservations. The demand elasticity has decayed to its floor. The population has adapted---not to the virus, but to the disruption.

Not rescue. Not collapse. A chronic condition the system learns to metabolise at declining cost, the way a household absorbs a family member's long illness: reorganising around the problem, never solving it, eventually barely noticing it. The pandemic becomes weather.

\needspace{4\baselineskip}
\bigskip

\noindent\textbf{Particle C: Managed Grind} \emph{(15.3\% of particles)}

\smallskip
\noindent Nearly four strains on average. The first wave peaks at 9.8\%, and as the waves accumulate, rejection begins its steady climb to 0.32. Not yet past the tipping point, but the vaccine programme is losing ground. The output gap troughs at $-0.95\%$. Mortality: 11.1\%.

Compare with B: similar strain count but the epidemic dynamics push harder. B accumulates 22.8 effective vaccine-weeks; C manages 15.7. The deeper troughs fire $f_5$ more forcefully, and the higher rejection thins the eligible arms at the same time. The virus outruns the platform from one direction. The population begins to drift from the other.

The output gap ends at $-0.66\%$, marginally worse than B's $-0.59\%$. The GDP trajectories are close. The public health trajectories are beginning to diverge. One in three refusing, versus one in three---but the trend line matters more than the terminal value.

\needspace{4\baselineskip}
\bigskip

\noindent\textbf{Particle D: Slow Ratchet} \emph{(21.0\% of particles)}

\smallskip
\noindent The single most probable cluster---one in five trajectories land here. Over four strains on average. The rejection ratchet has crossed the threshold: 0.38, nearly two in five refusing. Each wave adds a few more to the refusal column, and the vaccine programme thins to 10.6 effective weeks. The output gap troughs at $-0.81\%$.

This is the tipping point. Below this rejection level, the vaccine programme maintains enough coverage to blunt subsequent waves. Above it, the ratchet is self-reinforcing: more rejection means less vaccination, which means more infection, which means more mandates, which means more rejection. D sits on the boundary, the last cluster where the system is still fighting rather than settling.

Mortality: 12.5\%. Terminal output gap: $-0.70\%$, versus C's $-0.66\%$. The GDP difference is small. The structural difference is not: D has crossed into a regime where the vaccine programme is contracting rather than merely underperforming.

\needspace{4\baselineskip}
\bigskip

\noindent\textbf{Particle E: Structural Failure} \emph{(46.1\% of particles)}

\smallskip
\noindent Over six strains in three years on average. A new variant every five to six months, each one partially escaping the immunity built by its predecessor. The vaccine platform updates, dutifully, but each update is obsolete before the next shipment. The virus is not winning because it is more dangerous. It is winning because it is faster.

Rejection reaches 0.52. More than half refuse. The programme's effective coverage has collapsed: 7.2 effective vaccine-weeks over three years, against Particle~A's 23.3. The ratio is more than three to one.

Peak infection: 9.8\%. Cumulative mortality: 15.8\%---nearly double Particle~A's. The output gap troughs at $-1.14\%$ and ends at $-0.95\%$. Habituation prevents the $-3\%$ that constant coupling would produce, but the floor it finds is deep recession, not recovery. The economy is not collapsing. It has settled.

Nearly half of all simulated trajectories end here---the dominant outcome. The vaccine programme does not fail because it was badly designed. It fails because the conditions for its success---moderate rejection, economic stability, manageable strain evolution---were never simultaneously met. Each depends on the others. The coupling means that when one goes, it takes the rest. Habituation softens the blow but cannot reverse the direction. The compositional bearish bias, in this particle, is not an analytical artefact. It is the story.

\needspace{4\baselineskip}
\section{Four-Narrative Archetypal Trajectories}\label{app:narratives-4n}

Seven archetypal trajectories produced by $k$-medoids clustering on a 13-dimensional trajectory feature space (Section~\ref{sec:fiscal-results}), with each archetype represented by the weighted average of all particle paths in its cluster. Sorted by final rejection.

\needspace{4\baselineskip}
\bigskip

\noindent\textbf{Particle A: Mild Recovery} \emph{(13.3\% of particles)}

\smallskip
\noindent The world as the optimists imagine it. Roughly three strains in three years. Peak infection barely reaches 9.5\%, the output gap troughs at a shallow $-0.32\%$, and by the horizon it has clawed back to $-0.32\%$---close enough to zero that the pandemic, in retrospect, looks like a bad quarter. Rejection stays at 0.19, essentially unchanged. Nobody is marching against mandates because the mandates were barely needed. 25.4 effective vaccine-weeks---the highest of any particle by a wide margin. The vaccine works because nothing stops it from working.

The fiscal channel is a footnote. Peak spending reaches 8.0\% of GDP. The political ratchet does engage---activism rises to $\phi = 0.76$---but the recession is too shallow to translate high activism into large spending impulses, and the programme sunset ($g_{\text{decay}}$) erodes whatever spending accumulates. Terminal debt: 11.5~pp. The economy would have recovered without the spending. The government helped, modestly, with a problem that was already solving itself. This is what fiscal success looks like when the underlying conditions are benign---and it tells you nothing about what happens when they are not.

\needspace{4\baselineskip}
\bigskip

\noindent\textbf{Particle B: Fiscal Overshoot} \emph{(13.9\% of particles)}

\smallskip
\noindent The US particle. This is what actually happened---or close enough.

Five strains on average, a severe epidemic (9.4\% peak infection, 13.2\% mortality), and a fiscal response at full CARES Act scale: peak spending 8.4\% of GDP. The economy barely dips---minimum output gap $-0.08\%$, a recession so brief it would not show up in quarterly data---and then overshoots to $+2.66\%$ by the horizon. Positive. The only particle in the entire simulation, three-narrative or four-narrative, where the output gap ends above zero.

The overshoot is not an accident of mild strains. The epidemic is as severe as Particles~D or~E. What distinguishes~B is the combination of low rejection ($\rho_T = 0.19$, keeping the vaccine programme functional at 13.6 effective weeks) and massive, unimpeded fiscal activation. Low rejection means $f_{10}$ never fires---spending reaches the economy at full strength. The multiplier operates without friction.

And here is the punchline: the debt \emph{deviation} goes negative ($d_T = -47.9$~pp). A caveat on units: $d$ tracks deviation from the no-pandemic baseline, which was itself deteriorating---so a negative $d_T$ means the overshoot generates enough additional tax revenue to improve the fiscal position relative to where it was heading, not that the government is running a surplus. But the improvement is real, and it only happens because the intervention works. In every other particle, the same spending accumulates as additional debt.

Political activation stays at $\phi = 0.05$---the ratchet never engages because the recession never lasts long enough to trigger it. Self-negating intervention: the rescue succeeds precisely because it prevents the conditions that would have sustained it.

\needspace{4\baselineskip}
\bigskip

\noindent\textbf{Particle C: Delayed Grind} \emph{(6.6\% of particles)}

\smallskip
\noindent The particle that looks like B on paper and plays out like D.

Roughly three strains, peak infection at 9.4\%, fiscal spending reaching 7.9\% of GDP. The economy does not overshoot. The output gap troughs at $-0.53\%$ and grinds to $-0.50\%$.

Why? Timing and rejection. The waves land before the fiscal channel fully activates, and rejection reaches 0.21---moderate, but enough that $f_{10}$ begins to shave the fiscal multiplier. The economy absorbs the fiscal support and converts it into debt (18.1~pp) rather than growth. Same medicine, different patient, different outcome. The lowest-weight particle in the simulation: the system does not like knife edges.

\needspace{4\baselineskip}
\bigskip

\noindent\textbf{Particle D: Managed Decline} \emph{(30.2\% of particles)}

\smallskip
\noindent The modal particle. Nearly one in three trajectories lands here---the single most likely outcome in the four-narrative system. This is what ``normal'' looks like when normal includes a pandemic and a fiscal response.

Over four strains on average, moderate severity (9.5\% peak infection, 11.4\% mortality), rejection at 0.30. The economy contracts by 0.64\% at the trough and settles at $-0.64\%$. Peak spending: 8.0\% of GDP---not because political will is absent, but because the recession sits in a dead zone: deep enough to sustain spending, not catastrophic enough to trigger the emergency channel at full force. Fourteen effective vaccine-weeks. The programme functions, but at two-thirds capacity.

Debt: 20.6~pp. The fiscal channel buys a fraction of a percentage point of output gap---without it, this particle would be at $y_T \approx -0.9\%$---and the price is 21~pp of debt. Not a rescue. Not a failure. A sustained draw on the future for a partial present benefit. The managed decline: spending slows the fall and the bill arrives later. The adjacent particles---E and F, together another 31\%---end in the same neighbourhood, differing in degree but not in kind. Between them, D, E, and F account for more than 60\% of the distribution. The median trajectory ends in mild recession and moderate debt and the kind of exhausted stability that does not make for headlines but does make for a generation's lived experience.

\needspace{4\baselineskip}
\bigskip

\noindent\textbf{Particle E: Rejection Drag} \emph{(16.9\% of particles)}

\smallskip
\noindent Six strains on average, rejection at 0.39, mortality 16.6\%. The output gap troughs at $-0.75\%$ and ends at $-0.75\%$. Peak spending 7.9\% of GDP. The fiscal effort is there---the government is spending. It is spending into a headwind. $f_1$ drags the economy down through infection. $f_2$ disrupts supply. $f_5$ converts the recession into mandate backlash. One positive channel against three negative ones. The positive channel loses.

Debt: 23.1~pp. Effective vaccine-weeks: 8.5. The spending accumulates without producing recovery. Not collapse---the system is between basins, above the catastrophe of~G, below the managed mediocrity of~D. Enough refusal to impede recovery, not enough to break the system. The fiscal channel is running. It is running in place.

\needspace{4\baselineskip}
\bigskip

\noindent\textbf{Particle F: Fiscal Paralysis} \emph{(14.2\% of particles)}

\smallskip
\noindent Over six strains on average, rejection at 0.40, terminal output gap $-0.82\%$. The name is the diagnosis.

Peak spending: 9.5\% of GDP, the highest of any non-overshooting particle. The political ratchet is saturated ($\phi \approx 0.93$)---the system is \emph{trying} to spend. But rejection at 0.40 means $f_{10}$ reduces fiscal effectiveness by roughly 12\%. The money goes out. Twelve cents of every dollar are absorbed by the coupling structure before they reach the economy. The channel is not blocked. It is throttled.

Effective vaccine-weeks: 7.6. Debt: 28.3~pp. Without the fiscal channel, this particle would be at $y_T \approx -1.1\%$. The spending buys a quarter of a percentage point. At 28~pp of debt, the price per basis point of output gap is the highest of any particle that is not in outright structural failure. The fiscal channel works. It works expensively. And the population whose rejection is throttling it will bear the debt it accumulates.

\needspace{4\baselineskip}
\bigskip

\noindent\textbf{Particle G: Structural Failure} \emph{(5.0\% of particles)}

\smallskip
\noindent Over eight strains on average. Rejection at 0.57---more than one in two refusing. Mortality: 19.7\%. One in five dead. The output gap troughs at $-1.08\%$ and ends at $-1.05\%$. No recovery. The line is flat.

Peak spending: 8.3\% of GDP. With rejection above 0.5, $f_{10}$ blocks a sixth of every fiscal dollar. The emergency channel fires at full strength and it is not enough. The combined negative coupling---infection dragging demand, infection disrupting supply, recession feeding rejection---absorbs the fiscal impulse before it reaches the real economy. Debt: 34.8~pp. Effective vaccine-weeks: 4.3---a sixth of Particle~A's.

The four-narrative system produces the same structural failure as the three-narrative Particle~E, plus 35~pp of debt and no additional recovery. The fiscal channel, under these conditions, is pure cost. The system spends, borrows, and receives nothing, because every positive dollar is absorbed by the coupling structure before it arrives. This is the particle where the compositional bearish bias is not a bias at all. It is a description.



\begin{thebibliography}{10}

\bibitem[Fong(2016)]{fong2016}
Fong, B. (2016). \emph{The Algebra of Open and Interconnected Systems}. PhD thesis, University of Oxford.

\bibitem[Fong and Spivak(2019)]{fong2019}
Fong, B. and Spivak, D.~I. (2019). \emph{An Invitation to Applied Category Theory: Seven Sketches in Compositionality}. Cambridge University Press.

\bibitem[Gal\'i(2011)]{gali2011}
Gal\'i, J. (2011). Unemployment fluctuations and stabilization policies: a New Keynesian perspective. \emph{Zeuthen Lectures}. MIT Press.

\bibitem[Gal\'i(2015)]{gali2015}
Gal\'i, J. (2015). \emph{Monetary Policy, Inflation, and the Business Cycle}. Princeton University Press. 2nd edition.

\bibitem[Kermack and McKendrick(1927)]{kermack1927}
Kermack, W.~O. and McKendrick, A.~G. (1927). A contribution to the mathematical theory of epidemics. \emph{Proceedings of the Royal Society A}, 115(772):700--721.

\bibitem[Kschischang et~al.(2001)]{kschischang2001}
Kschischang, F.~R., Frey, B.~J., and Loeliger, H.-A. (2001). Factor graphs and the sum-product algorithm. \emph{IEEE Transactions on Information Theory}, 47(2):498--519.

\bibitem[Phoa(2026)]{phoa2026}
Phoa, W. (2026). Particle filters and factor graphs for narrative space: a compositional framework for scenario analysis and asset allocation. In preparation.

\bibitem[Taylor(1993)]{taylor1993}
Taylor, J.~B. (1993). Discretion versus policy rules in practice. \emph{Carnegie-Rochester Conference Series on Public Policy}, 39:195--214.

\end{thebibliography}
\end{document}